\documentclass[12pt]{amsart}
\usepackage{amscd}      
\usepackage{amssymb}
\usepackage{amsmath}
\usepackage{xypic}      
\LaTeXdiagrams          
\usepackage[all,v2]{xy}
\xyoption{2cell} \UseAllTwocells \xyoption{frame} \CompileMatrices
\allowdisplaybreaks[3]

\usepackage{latexsym}
\usepackage{epsfig}
\usepackage{amsfonts}
\usepackage{enumerate}

\newtheorem{prop}{Proposition}[section]
\newtheorem{lem}[prop]{Lemma}

\newtheorem{cor}[prop]{Corollary}
\newtheorem{thm}[prop]{Theorem}

\numberwithin{equation}{section}

\newtheorem{defn}[prop]{Definition}

\newtheorem{example}[prop]{Example}

\newtheorem{rmk}[prop]{Remark}

\newenvironment{pf}{\begin{trivlist}\item[]{\sc Proof.}}%
            {\nolinebreak $\Box$ \end{trivlist}}

\newcommand{\noprint}[1]{}

\newcommand{\toto}{\rightrightarrows}

\newcommand{\ldiag}[1]%
       {\makebox[0cm]{${\scriptstyle#1}\downarrow\phantom{\scriptstyle#1}$}}
\newcommand{\ldiagup}[1]%
       {\makebox[0cm]{${\scriptstyle#1}\uparrow\phantom{\scriptstyle#1}$}}
\newcommand{\rdiag}[1]%
       {\makebox[0cm]{$\phantom{\scriptstyle#1}\downarrow{\scriptstyle#1}$}}
\newcommand{\sediagr}[1]%
       {\makebox[0cm]{$\phantom{\scriptstyle#1}\searrow{\scriptstyle#1}$}}
\newcommand{\nediagr}[1]%
       {\makebox[0cm]{$\phantom{\scriptstyle#1}\nearrow{\scriptstyle#1}$}}
\newcommand{\rdiagup}[1]%
       {\makebox[0cm]{$\phantom{\scriptstyle#1}\uparrow{\scriptstyle#1}$}}
\newcommand{\swdiag}[1]%
       {\makebox[0cm]{$\phantom{\scriptstyle#1}\swarrow{\scriptstyle#1}$}}
\newcommand{\sediag}[1]%
       {\makebox[0cm]{${\scriptstyle#1}\searrow\phantom{\scriptstyle#1}$}}
\newcommand{\nediag}[1]%
       {\makebox[0cm]{${\scriptstyle#1}\nearrow\phantom{\scriptstyle#1}$}}

\newcommand{\doublearrowstack}[2]%
                      {{{{\scriptstyle#1}\atop{\textstyle\longrightarrow}}\atop{{\textstyle\longrightarrow}\atop{\scriptstyle#2}}}}
\newcommand{\rightleftarrowstack}[2]%
                      {{{{\scriptstyle#1}\atop{\textstyle\longrightarrow}}\atop{{\textstyle\longleftarrow}\atop{\scriptstyle#2}}}}
\newcommand{\leftrightarrowstack}[2]%
                      {{{{\scriptstyle#1}\atop{\textstyle\longleftarrow}}\atop{{\textstyle\longrightarrow}\atop{\scriptstyle#2}}}}

\newcommand{\overtoparrow}%
{\makebox[0cm]{\beginpicture \setcoordinatesystem units
<.8cm,.4cm> point at 0 0 \setplotarea x from -3 to 3, y from 0 to
1 \setquadratic \plot -3 0 0 1 3 0 / \put{\vector(3,-1){0}}[Bl] at
3 0
\endpicture}}

\newcommand{\underbottomarrow}%
{\makebox[0cm]{\beginpicture \setcoordinatesystem units
<.8cm,.4cm> point at 0 0 \setplotarea x from -3 to 3, y from 0 to
1 \setquadratic \plot -3 1 0 0 3 1 / \put{\vector(3,1){0}}[Bl] at
3 1
\endpicture}}

\newcommand{\ses}[5]%
{0\longrightarrow#1\stackrel{#2}{ \longrightarrow}#3\stackrel{#4}{
\longrightarrow}#5\longrightarrow0}

\newcommand{\dt}[6]%
{#1\stackrel{#2}{longrightarrow}#3
\stackrel{#4}{\longrightarrow}#5 \stackrel{#6}{\longrightarrow}
#1[1]}

\newcommand{\cat}[1]%
{(\mbox{\rm #1})}

\def\Label#1{\label{#1}{\tt [#1]}\phantom{h}}
\def\Label{\label}

\title[Orbifold Chow rings for hypertoric DM stacks]{The Orbifold Chow Ring of Hypertoric Deligne-Mumford Stacks}

\author{Yunfeng Jiang}
\address{Department of Mathematics\\ University of British Columbia\\ 1984 Mathematics Road\\Vancouver\\ BC V6T 1Z2\\ Canada}
\email{jiangyf@math.ubc.ca}
\author{Hsian-Hua Tseng}
\address{Department of Mathematics\\ University of British Columbia\\ 1984 Mathematics Road\\Vancouver\\ BC V6T 1Z2\\ Canada}
\email{hhtseng@math.ubc.ca}
\date{\today}

\begin{document}
\begin{abstract}
Hypertoric varieties are determined by hyperplane arrangements. In this paper,
we use stacky hyperplane arrangements to define the notion of hypertoric Deligne-Mumford stacks.
Their orbifold Chow rings are computed.
As an application, some examples related to crepant resolutions are discussed.
\end{abstract}

\maketitle

\section{Introduction}
Hypertoric varieties (cf. \cite{BD}, \cite{p}) are the hyperk\"ahler analogue of K\"ahler
toric varieties. The algebraic
construction of hypertoric varieties was given by Hausel and
Sturmfels \cite{HS}. Modelling on their construction, in this paper
we construct hypertoric DM stacks and study their
orbifold Chow rings.

According to \cite{BD}, the topology of hypertoric varieties is determined by hyperplane arrangements. 
In this paper we define stacky hyperplane arrangements from which we 
define hypertoric DM stacks.

Let $N$ be a finitely generated abelian group of rank $d$ and
$N\to \overline{N}$ the natural projection modulo
torsion. Let $\beta:
\mathbb{Z}^{m}\to N$ be a homomorphism determined by a
collection of nontorsion integral vectors $\{b_{1},\cdots,b_{m}\}\subseteq
N$. We require that $\beta$ has finite cokernel. The Gale dual of $\beta$ is denoted by $\beta^{\vee}:
(\mathbb{Z}^{m})^{*}\to DG(\beta)$. A $generic$ element 
$\theta$ in $DG(\beta)$ and the vectors $\{\overline{b}_{1},\cdots,\overline{b}_{m}\}$
determine a hyperplane arrangement $\mathcal{H}=(H_{1},\cdots,H_{m})$ in
$N_{\mathbb{R}}^{*}$. We call $\mathcal{A}:=(N,\beta,\theta)$
a {\em stacky hyperplane arrangement}.

For $\beta:
\mathbb{Z}^{m}\to N$  in $\mathcal{A}$, we consider the Lawrence lifting
$\beta_{L}: \mathbb{Z}^{m}\oplus \mathbb{Z}^{m}\rightarrow N_{L}$ of $\beta$
where $N_{L}$ is a finitely generated abelian group with rank $m+d$. The map $\beta_{L}$ is 
given by vectors $\{b_{L,1},\cdots,b_{L,m},b^{'}_{L,1},\cdots,b^{'}_{L,m}\}\subseteq
N_{L}$. The generic element $\theta$ determines a Lawrence simplicial fan $\Sigma_{\theta}$
in $\overline{N}_{L}$. We call $\mathbf{\Sigma_{\theta}}=(N_{L},\Sigma_{\theta},\beta_{L})$ a Lawrence stacky fan and
$\mathcal{X}(\mathbf{\Sigma_{\theta}})$ the Lawrence toric DM stack.
The hypertoric DM stack  $\mathcal{M}(\mathcal{A})$
associated to $\mathcal{A}$ is defined as  a quotient stack which is a closed substack of the Lawrence toric DM stack 
$\mathcal{X}(\mathbf{\Sigma_{\theta}})$, generalizing 
the construction of \cite{HS}. 

The stacky hyperplane arrangement $\mathcal{A}$ also  determines an extended 
stacky fan $\mathbf{\Sigma}=(N,\Sigma,\beta)$  
introduced in \cite{Jiang}. Here $\Sigma$ is the normal fan of the bounded
polytope $\mathbf{\Gamma}$ of the hyperplane arrangement $\mathcal{H}$.
The toric DM stack $\mathcal{X}(\mathbf{\Sigma})$ defined in \cite{Jiang} is 
the associated toric DM stack of $\mathcal{M}(\mathcal{A})$.

To the map $\beta $  we associate a
multi-fan $\Delta_{\mathbf{\beta}}$ in the sense of \cite{HM}, which consists of  cones generated by linearly independent subsets
$\{\overline{b}_{i_{1}},\cdots,\overline{b}_{i_{k}}\}$ in $\overline{N}$ for 
$\{i_{1},\cdots,i_{k}\}\subset\{1,\cdots,m\}$,
see Section \ref{substack}. We assume that the $supp(\Delta_{\mathbf{\beta}})=\overline{N}$. 
We prove that each top dimensional cone in 
$\Delta_{\mathbf{\beta}}$ gives a local chart for the hypertoric DM stack 
$\mathcal{M}(\mathcal{A})$. We define a set $Box(\Delta_{\mathbf{\beta}})$ consisting of all 
pairs $(v,\sigma)$, where $\sigma$ is a cone in the multi-fan
$\Delta_{\mathbf{\beta}}$,  $v\in N$ such that 
$\overline{v}=\sum_{\rho_{i}\subset\sigma}\alpha_{i}\overline{b}_{i}$ for 
$0<\alpha_{i}<1$.  For $(v,\sigma)\in
Box(\Delta_{\mathbf{\beta}})$ we consider a closed substack of
$\mathcal{M}(\mathcal{A})$ given by the quotient stacky hyperplane arrangement
$\mathcal{A}(\sigma)$. The inertia stack of  $\mathcal{M}(\mathcal{A})$ is
the disjoint union of all such closed substacks, see Section
\ref{substack}.

We now describe the orbifold Chow ring of
$\mathcal{M}(\mathcal{A})$. The multi-fan $\Delta_{\mathbf{\beta}}$ naturally gives a ``matroid''
$M_{\beta}$. The vertex set is $\{1,\cdots,m\}$, and the faces are the subsets 
$\{i_1,\cdots,i_k\}\subseteq \{1,\cdots,m\}$ such that $\{\overline{b}_{i_{1}},\cdots,\overline{b}_{i_{k}}\}$
are linearly independent in $\overline{N}$. Note that the faces of $M_{\beta}$ are the  cones in $\Delta_{\mathbf{\beta}}$. 
According to \cite{HS}, the ordinary cohomology ring of the hypertoric variety corresponding to the hyperplane arrangement $\mathcal{H}$ is isomorphic to the  ``Stanley-Reisner'' ring of the matroid $M_{\beta}$. Our result shows that the orbifold Chow ring of hypertoric DM stacks is a generalization of 
the  Stanley-Reisner ring of the matroid $M_{\beta}$ to the multi-fan $\Delta_{\mathbf{\beta}}$. 
Let $N^{\Delta_\mathbf{\beta}}$ denote  
all the pairs $(c,\sigma)$, where $c\in N$,  $\sigma$ is a cone in  
$\Delta_\beta$ such that 
$\overline{c}=\sum_{\rho_{i}\subseteq \sigma}a_{i}\overline{b}_{i}$ and  
$a_{i}>0$ are  rational numbers.  Then $N^{\Delta_\mathbf{\beta}}$ gives rise to a 
group ring
$$\mathbb{Q}[\Delta_\mathbf{\beta}]=\bigoplus_{(c,\sigma)\in N^{\Delta_\mathbf{\beta}}}\mathbb{Q}\cdot y^{(c,\sigma)},$$
where $y$ is a formal variable.  For any $(c,\sigma)\in N^{\Delta_\mathbf{\beta}}$, there exists a unique element
$(v,\tau)\in Box(\Delta_\mathbf{\beta})$ such that $\tau\subset\sigma$ and 
$c=v+\sum_{\rho_{i}\subseteq \sigma}m_{i}b_{i}$, 
where 
$m_{i}$ are nonnegative integers. We call $(v,\tau)$ the $fractional~ part$ of $(c,\sigma)$. For $(c,\sigma)$ we define the {\em ceiling function} $\lceil c \rceil_{\sigma}$ by  
$\lceil c \rceil_{\sigma}=\sum_{\rho_{i}\subseteq \tau}b_{i}+\sum_{\rho_{i}\subseteq \sigma}m_{i}b_{i}$. Note that 
if  $\overline{v}=0$, $\lceil c \rceil_{\sigma}=\sum_{\rho_{i}\subseteq \sigma}m_{i}b_{i}$.   
For two pairs  $(c_1,\sigma_1)$, $(c_2,\sigma_2)$, if $\sigma_{1}\cup\sigma_{2}$ is a cone in $\Delta_\mathbf{\beta}$, define 
$\epsilon(c_1,c_2):=\lceil c_1 \rceil_{\sigma_{1}}+\lceil c_2 \rceil_{\sigma_{2}}-\lceil c_1+c_2 \rceil_{\sigma_{1}\cup\sigma_2}$.
Let $\sigma_{\epsilon}\subseteq\sigma_1\cup\sigma_2$ be the minimal cone in $\Delta_\mathbf{\beta}$ containing $\epsilon(c_1,c_2)$ so that 
$(\epsilon(c_1,c_2),\sigma_{\epsilon})\in N^{\Delta_\mathbf{\beta}}$. We define the grading on $\mathbb{Q}[\Delta_\mathbf{\beta}]$ as follows.
For any  $(c,\sigma)$, write $c=v+\sum_{\rho_{i}\subseteq \sigma}m_{i}b_{i}$,  then
$$deg(y^{(c,\sigma)}):=|\tau|+\sum_{\rho_{i}\subseteq\sigma}m_{i},$$ where $|\tau|$ is the dimension of $\tau$. 
Note that we have integer grading due to the fact that $\mathcal{M}(\mathcal{A})$ is hyperk\"ahler. By abuse of notation, we write $y^{(b_{i},\rho_i)}$ as $y^{b_{i}}$.
The multiplication 
is defined by 
\begin{equation}\Label{product}
y^{(c_{1},\sigma_{1})}\cdot y^{(c_{2},\sigma_{2})}:=
\begin{cases}
(-1)^{|\sigma_{\epsilon}|}y^{(c_{1}+c_{2}+\epsilon(c_1,c_2),\sigma_{1}\cup\sigma_{2})}&\text{if
$\sigma_{1}\cup\sigma_{2}$ is a cone in $\Delta_{\mathbf{\beta}}$}\,,\\
0&\text{otherwise}\,.
\end{cases} 
\end{equation}
Using the property of ceiling functions we check that the multiplication is commutative and associative. So $\mathbb{Q}[\Delta_\mathbf{\beta}]$ is 
a unital associative commutative	ring. 
Let
$Cir(\Delta_\mathbf{\beta})$ be the ideal in
$\mathbb{Q}[\Delta_\mathbf{\beta}]$ generated by the elements:
\begin{equation}\Label{ideal2}
\sum_{i=1}^{m}e(b_{i})y^{b_{i}}, \quad e\in
N^{*}. 
\end{equation} 
Let
$A^{*}_{orb}(\mathcal{M}(\mathcal{A}))$ be the orbifold Chow
ring of the hypertoric DM stack $\mathcal{M}(\mathcal{A})$. We
have the following Theorem:

\begin{thm}\label{chowring}
Let $\mathcal{M}(\mathcal{A})$ be the hypertoric DM stack associated to the  stacky hyperplane arrangement $\mathcal{A}$.
Then there is an isomorphism  of graded $\mathbb{Q}$-algebras:
$$A^{*}_{orb}(\mathcal{M}(\mathcal{A}))\cong \frac{\mathbb{Q}[\Delta_\mathbf{\beta}]}{Cir(\Delta_\mathbf{\beta})}.$$
\end{thm}
The orbifold Chow ring of the hypertoric DM stack $\mathcal{M}(\mathcal{A})$ is independent 
of the generic element $\theta$. It only depends on the map $\beta$.

Theorem \ref{chowring} is proven by a direct approach. The inertia
stack of a hypertoric DM stack $\mathcal{M}(\mathcal{A})$ is
the disjoint union of closed substacks
$\mathcal{M}(\mathcal{A}(\sigma))$ for all $(v,\sigma)\in
Box(\Delta_{\mathbf{\beta}})$. To determine the ring structure,
we identify the 3-twisted sectors as closed substacks of
$\mathcal{M}(\mathcal{A})$ indexed by triples
$((v_{1},\sigma_{1}),(v_{2},\sigma_{2}),(v_{3},\sigma_{3}))$ in $Box(\Delta_{\mathbf{\beta}})^3$ such that 
$v_{1}+v_{2}+v_{3}\in N$ is a integral linear combination of $b_i$'s. We then determine the obstruction bundle
over any 3-twisted sector and prove that the orbifold cup product is the same as the product of the ring
$\mathbb{Q}[\Delta_\mathbf{\beta}]$ described above.

The multi-fan $\Delta_\mathbf{\beta}$ is equal to the simplicial fan $\Sigma$ in $\mathbf{\Sigma}$ induced 
from the stacky hyperplane arrangement $\mathcal{A}$ if and only if $\mathcal{H}$
has $n$ hyperplanes $\{H_{1},\cdots,H_{n}\}$ whose normal polytope is a product of simplices. So in this case 
$\mathbf{\Sigma}$ is a stacky fan and the simplicial fan $\Sigma$ is a product of normal
fans of simplices, the toric variety $X(\Sigma)$ is a product of weighted projective
spaces. Then by \cite{BD} the associated hypertoric variety is the
cotangent bundle of the toric variety $X(\Sigma)$. So
$\mathcal{M}(\mathcal{A})\simeq
T^*\mathcal{X}(\mathbf{\Sigma})$, the cotangent bundle of the toric DM stack
$\mathcal{X}(\mathbf{\Sigma})$.  The ring 
$\mathbb{Q}[\Delta_\mathbf{\beta}]$ coincides (as vector spaces) with the deformed ring
$\mathbb{Q}[N]^\mathbf{\Sigma}$ as defined in
\cite{BCS}.  

\begin{cor}
Let $\mathbf{\Sigma}$ be as above. Then there is an isomorphism of
$\mathbb{Q}$-vector spaces
$$A_{orb}^*(\mathcal{M}(\mathcal{A}))\simeq
A_{orb}^*(\mathcal{X}(\mathbf{\Sigma})).$$
\end{cor}

Here is an example which shows that the  orbifold Chow ring of
$\mathcal{M}(\mathcal{A})$ is not isomorphic as a ring to
the orbifold Chow ring of the associated toric DM stack
$\mathcal{X}(\mathbf{\Sigma})$. Consider the weighted projective stack $\mathbb{P}(1,2)$ which is a toric DM stack with stacky fan 
$\mathbf{\Sigma}=(N,\Sigma,\beta)$, where $N=\mathbb{Z}$, $\beta: \mathbb{Z}^{2}\rightarrow N$ is given by the 
vectors $b_1=(1),b_2=(-2)$ and $\Sigma$ is the simplicial fan in the lattice $N$ consisting cones $\rho_1$ and $\rho_2$ generated by 
$b_1=(1)$ and $b_2=(-2)$ respectively. The Gale dual map $\beta^{\vee}: \mathbb{Z}^{2}\rightarrow \mathbb{Z}$ is given by 
the matrix $(2)$. Choosing generic element $\theta=(1)\in \mathbb{Z}$, we get a stacky hyperplane arrangement 
$\mathcal{A}=(N,\beta,\theta)$. The hypertoric DM stack $\mathcal{M}(\mathcal{A})$ is the  cotangent bundle $T^{*}\mathbb{P}(1,2)$  whose core is the toric DM stack $\mathbb{P}(1,2)$. Both $\mathbb{Q}[\Delta_{\beta}]$ and $\mathbb{Q}[N]^{\mathbf{\Sigma}}$ are generated by 
$y^{b_1}$, $y^{b_2}$, and $y^{(\frac{1}{2}b_2,\rho_2)}$. According to Theorem 1.1 and the main theorem in \cite{BCS}, their orbifold Chow rings are given as follows:
$$
A^{*}_{orb}(\mathcal{X}(\mathbf{\Sigma});\mathbb{Q})\cong \frac{\mathbb{Q}[x_{1},x_{2},v]}
{(x_{1}-2x_{2},v^{2}-x_{2},vx_{1},x_{1}x_{2})}
\cong \frac{\mathbb{Q}[v]}
{(v^{3})},
$$
$$
A^{*}_{orb}(\mathcal{M}(\mathcal{A});\mathbb{Q})\cong \frac{\mathbb{Q}[x_{1},x_{2},v]}
{(x_{1}-2x_{2},x_{1}x_{2},vx_{1},v^{2})}
\cong \frac{\mathbb{Q}[x_{2},v]}
{(x_{2}^{2},vx_{2},v^{2})}.
$$
It is easy to see that these two rings are not isomorphic. Thus the orbifold Chow ring
of a hypertoric DM stack is not necessarily isomorphic to the orbifold Chow ring 
of its core. (However, their Chow rings are isomorphic, see Theorem 1.1 of \cite{HS}.) 
This also proves that the orbifold Chow ring has no homotopy invariance property. We remark 
that the core of a general hypertoric DM stack can be singular, it is not clear how to 
define orbifold Chow ring. But in the case of a cotangent bundle over weighted projective space,
the core is the weighted projective space and the orbifold Chow ring is well-defined.
On the other hand, the orbifold Chow ring of a Lawrence toric DM stack is isomorphic to 
its associated hypertoric DM stack, see \cite{JT}.

Computations of orbifold cohomology rings of hypertoric orbifolds
in symplectic geometry have been pursued in \cite{GH}.

This paper is organized as follows. In Section \ref{hyperplane} we
discuss the relation  between  stacky hyperplane arrangements and
extended stacky fans. We define
hypertoric DM stack $\mathcal{M}(\mathcal{A})$ associated to
the stacky hyperplane arrangement $\mathcal{A}$. In Section \ref{hypertoricDM}
we discuss the properties of hypertoric DM stacks.  In Section
\ref{substack} we determine closed substacks of a hypertoric DM
stack. This yields a description of its inertia stacks. We prove
Theorem \ref{chowring} in Section \ref{chow}, and in Section
\ref{app} we give some examples.
\subsection*{Conventions}
In this paper we work entirely algebraically over the field of
complex numbers. Chow rings and orbifold Chow rings are taken with
rational coefficients. By an orbifold we mean a smooth
Deligne-Mumford stack with trivial generic stabilizer.
We refer to \cite{BCS} for the construction of Gale dual
$\beta^{\vee}: \mathbb{Z}^{m}\to DG(\beta)$ from 
$\beta: \mathbb{Z}^{m}\to N$. We denote by $N\to \overline{N}$ the natural map
modulo torsion.  For cones $\sigma_{1},\sigma_{2}$ in $\mathbb{R}^{d}$,
we use $\sigma_{1}\cup\sigma_{2}$ to represent the set of union of the 
generators of $\sigma_{1}$ and $\sigma_{2}$.

\subsection*{Acknowledgments}
We would like to thank the referee for the nice comments of Remark \ref{lawrencermk-s}
and Kai Behrend, Megumi Harada, Nicholas Proudfoot for valuable discussions.
\section{The Hypertoric DM Stacks}\label{hyperplane}
In this section we define  hypertoric Deligne-Mumford stacks,
mimicking the construction of hypertoric varieties in \cite{HS}.
\subsection*{Stacky hyperplane arrangements}
We introduce stacky hyperplane arrangements.  We explain how a stacky hyperplane arrangement
gives extended stacky fans. 

Let $N$ be a finitely generated abelian group and  $\beta:\mathbb{Z}^m\to N$  a map given by nontorsion
integral vectors $\{b_1,...,b_m\}$. We have the following exact
sequences:
\begin{equation}\Label{exact1}
0\longrightarrow DG(\beta)^{*}\stackrel{(\beta^{\vee})^{*}}{\longrightarrow}
\mathbb{Z}^{m}\stackrel{\beta}{\longrightarrow} N\longrightarrow
Coker(\beta)\longrightarrow 0,
\end{equation}
\begin{equation}\Label{exact2}
0\longrightarrow N^{*}\longrightarrow
\mathbb{Z}^{m}\stackrel{\beta^{\vee}}{\longrightarrow}
DG(\beta)\longrightarrow Coker(\beta^{\vee})\longrightarrow 0,
\end{equation}
where $\beta^{\vee}$ is the Gale dual of $\beta$ (see \cite{BCS}). 
The map $\beta^{\vee}$ is given by the integral vectors $\{a_1,\cdots,a_m\}\subseteq DG(\beta)$. 
Choose a generic element $\theta\in DG(\beta)$ which lies in the inage of $\beta^{\vee}$ and let
$\psi:=(r_{1},\cdots,r_{m})$ be a lifting of $\theta$ in
$\mathbb{Z}^{m}$ such that $\theta=-\beta^{\vee}\psi$. Note that $\theta$ is
generic if and only if it is not in any hyperplane of the
configuration determined by $\beta^{\vee}$ in $DG(\beta)_{\mathbb{R}}$. Let 
$M=N^{*}$ be the dual of $N$ and $M_{\mathbb{R}}=M\otimes_{\mathbb{Z}}\mathbb{R}$, then
$M_{\mathbb{R}}$ is a $d$-dimensional $\mathbb{R}$-vector space. 
Associated to $\theta$ there is a hyperplane arrangement
$\mathcal{H}=\{H_{1},\cdots,H_{m}\}$ in $M_{\mathbb{R}}$ defined by
$H_{i}$  the hyperplane
\begin{equation}\Label{arrangement}
H_{i}:=\{v\in M_{\mathbb{R}}|<b_{i},v>+r_{i}=0\}\subset M_{\mathbb{R}}.
\end{equation}
So (\ref{arrangement}) determines hyperplane arrangements in $M_{\mathbb{R}}$, up to translation
induced by the choice of the lifting $\psi:=(r_{1},\cdots,r_{m})$.

\begin{defn}
We call
$\mathcal{A}:=(N,\beta,\theta)$ a {\em stacky hyperplane arrangement}.
\end{defn}

It is
well-known that hyperplane arrangements determine the topology of
hypertoric varieties \cite{BD}. 
Let
$$\mathbf{\Gamma}=\bigcap_{i=1}^{m}F_{i}, \text{ where }F_{i}=\{v\in M_{\mathbb{R}}|<b_{i},v>+r_{i}\geq 0\}.$$ 
Let $\Sigma$ be the normal fan of 
$\mathbf{\Gamma}$ in $M_{\mathbb{R}}=\mathbb{R}^{d}$ with one
dimensional rays generated by
$\overline{b}_{1},\cdots,\overline{b}_{n}$. By reordering, we may
assume that $H_{1},\cdots,H_{n}$ are the hyperplanes that bound
the polytope $\mathbf{\Gamma}$, and $H_{n+1},\cdots,H_{m}$ are the
other hyperplanes. Then we have an extended stacky fan
$\mathbf{\Sigma}=(N,\Sigma,\beta)$ defined in \cite{Jiang}, where
$\beta:
\mathbb{Z}^{m}	\to N$ is given by
$\{b_{1},\cdots,b_{n},b_{n+1},\cdots,b_{m}\}\subset N$, and  $\{b_{n+1},\cdots,b_{m}\}$
are the extra data. 

By \cite{Jiang}, the extended stacky fan $\mathbf{\Sigma}$
determines a toric Deligne-Mumford stack
$\mathcal{X}(\mathbf{\Sigma})$. It is the same stack as in \cite{BCS}.
Its coarse moduli space is the
toric variety corresponding to the normal fan of
$\mathbf{\Gamma}$. According to \cite{BD}, a hyperplane arrangement
$\mathcal{H}$ is {\em simple} if the codimension of the nonempty
intersection of  any $l$ hyperplanes is $l$. A hypertoric variety
is the coarse moduli space of an  {\em orbifold} if the corresponding hyperplane arrangement
is simple.

\begin{example}\Label{ex}
Let $\mathcal{H}=\{H_{1},H_{2},H_{3},H_{4}\}$, see Figure 1. The polytope $\mathbf{\Gamma}$ of the hyperplane
arrangement is the shaded triangle whose toric variety is the projective
plane. The extended stacky fan is given by the fan of the
projective plane $\mathbb{P}^{2}$ and an extra ray $(0,1)$.
\begin{center}
\begin{picture}(0,0)%
\includegraphics{figure1.pstex}%
\end{picture}%
\setlength{\unitlength}{3947sp}%
\begingroup\makeatletter\ifx\SetFigFont\undefined%
\gdef\SetFigFont#1#2#3#4#5{%
  \reset@font\fontsize{#1}{#2pt}%
  \fontfamily{#3}\fontseries{#4}\fontshape{#5}%
  \selectfont}%
\fi\endgroup%
\begin{picture}(4007,2064)(3901,-3940)
\end{picture}

\end{center}
\end{example}

\begin{rmk}
If for a generic element $\theta\in DG(\beta)$ the
hyperplane arrangement $\mathcal{H}$ bounds a polytope whose
normal fan is $\Sigma$, then $\mathbf{\Sigma}=(N,\Sigma,\beta)$ is a stacky fan defined in
\cite{BCS}. 
\end{rmk}

\subsection*{Lawrence toric DM stacks}
Consider the Gale dual map
$\beta^{\vee}: \mathbb{Z}^{m}\to
DG(\beta)$ in (\ref{exact2}).  We denote 
the Gale dual map of 
$$ \mathbb{Z}^{m}\oplus \mathbb{Z}^{m}\stackrel{(\beta^{\vee},-\beta^{\vee})}{\longrightarrow}
DG(\beta)$$ 
by 
\begin{equation}\Label{betal}
\beta_{L}: \mathbb{Z}^{2m}\rightarrow N_{L},
\end{equation} \\
where $\overline{N}_{L}$
is a lattice of dimension $2m-(m-d)$. The map $\beta_{L}$ is given by the integral vectors  $\{b_{L,1},\cdots,
b_{L,m},b'_{L,1},\cdots,
b^{'}_{L,m}\}$ and $\beta_{L}$ is called  the Lawrence lifting of $\beta$.
\begin{rmk}\label{lawrencermk-s}
Consider the following commutative diagram for which the rows are exact:
$$\xymatrix{
~\mathbb{Z}^{m}\ar_{[1,1]^{t}}@<-.5ex>[r]\dto & ~\mathbb{Z}^{2m}\ar_{[0,1]}@<-.5ex>[l]\ar_{[1,-1]}@<-.5ex>[r]\dto^{(\beta^{\vee},-\beta^{\vee})} & ~\mathbb{Z}^{m}\ar_{[1,0]^{t}}
@<-.5ex>[l]\dto^{\beta^{\vee}}\\
0\ar@<-.5ex>[r] & DG(\beta)\ar@<-.5ex>[l]\ar@<-.5ex>[r]& DG(\beta),\ar@<-.5ex>[l]}
$$
where $[1,1]^{t}$ represents the transpose of the matrix. Taking Gale dual to the above diagram yields
the following commutative diagram:
$$\xymatrix{
~\mathbb{Z}^{m}\ar^{[0,1]^{t}}@<.5ex>[r]\dto^{id} & ~\mathbb{Z}^{2m}\ar^{[1,1]}@<.5ex>[l]\ar^{[1,0]}@<.5ex>[r]\dto^{\beta_{L}} & ~\mathbb{Z}^{m}\ar^{[1,-1]^{t}}@<.5ex>[l]\dto^{\beta}\\
~\mathbb{Z}^{m}\ar@<.5ex>[r] & N_{L}\ar@<.5ex>[l]\ar@<.5ex>[r]& N.\ar@<.5ex>[l]}
$$
So from the functoriality of Gale dual we have $N_{L}\cong N\oplus\mathbb{Z}^{m}$.
Since $(\beta^{\vee},-\beta^{\vee})=\beta^{\vee}\oplus 0$, we get that $\beta_{L}=\beta\oplus id$.
So $\{b_{L,1},\cdots,
b_{L,m},b'_{L,1},\cdots,
b^{'}_{L,m}\}$ are the vectors $$\{(b_{1},e_{1}),\cdots,
(b_{m},e_{m}),(0,e_{1}),\cdots,
(0,e_{m})\},$$ where $\{e_{i}\}$ are the standard bases of $\mathbb{Z}^{m}$.
\end{rmk}

Given the generic element $\theta$, let $\overline{\theta}$
be the natural image of $\theta$ under the projection $DG(\beta)\rightarrow \overline{DG(\beta)}$. 
Then the map $\overline{\beta}^{\vee}: \mathbb{Z}^{m}\rightarrow \overline{DG(\beta)}$ is given by $\overline{\beta}^{\vee}=(\overline{a}_{1},\cdots,\overline{a}_{m})$. For any basis of $\overline{DG(\beta)}$
of the form $C=\{\overline{a}_{i_{1}},\cdots,\overline{a}_{i_{m-d}}\}$, there exist unique
$\lambda_{1},\cdots,\lambda_{m-d}$ such that
$$\overline{a}_{i_{1}}\lambda_{1}+\cdots+\overline{a}_{i_{m-d}}\lambda_{m-d}=\overline{\theta}.$$
Let $\mathbb{C}[z_{1},\cdots,z_{m},w_{1},\cdots,w_{m}]$ be the coordinate ring of $\mathbb{C}^{2m}$. Let
$$\sigma(C,\theta)=\{\overline{b}_{L,i_{j}}~|\lambda_{j}>0\}\sqcup\{\overline{b}'_{L,i_{j}}|~\lambda_{j}<0\} \quad \text{and} \quad 
C(\theta)=\{z_{i_{j}}~|\lambda_{j}>0\}\sqcup\{w_{i_{j}}|~\lambda_{j}<0\}.$$
We put
\begin{equation}\Label{irrelevant}
\mathbf{\mathcal{I}}_{\theta}:=\left\langle\prod
C(\theta)|~C~\text{is a basis of}~\overline{DG(\beta)}\right\rangle,
\end{equation}
and 
\begin{equation}\Label{fan}
\Sigma_{\theta}:=\{\overline{\sigma}(C,\theta):~C~\text{is a basis of}~\overline{DG(\beta)}\},
\end{equation}
where $\overline{\sigma}(C,\theta)=\{\overline{b}_{L,1},\cdots,\overline{b}_{L,m},\overline{b}'_{L,1},\cdots,\overline{b}'_{L,m}\}\setminus\sigma(C,\theta)$ 
is the complement of $\sigma(C,\theta)$ and corresponds to a maximal cone in $\Sigma_{\theta}$.
From \cite{HS}, $\Sigma_{\theta}$ is the fan of a Lawrence toric variety $X(\Sigma_{\theta})$ corresponding to
$\theta$ in the lattice $\overline{N}_{L}$, and $\mathcal{I}_{\theta}$ is the irrelevant ideal.  
The construction above establishes the following
\begin{prop}\Label{lawrencestacky}
A stacky hyperplane arrangement $\mathcal{A}=(N,\beta,\theta)$
also gives a stacky fan $\mathbf{\Sigma_{\theta}}=(N_{L},\Sigma_{\theta},\beta_{L})$
which is called a Lawrence stacky fan.  
\end{prop}

\begin{pf}
From Proposition 4.3 in \cite{HS}, $\Sigma_{\theta}$ is a simplicial fan in $\overline{N}_{L}$.  The rays $\rho_{L,i}$, $\rho^{'}_{L,i}$
are generated by $\overline{b}_{L,i}$, $\overline{b}^{'}_{L,i}$. The map $\beta_{L}$ is the map 
(\ref{betal}) given by $\{b_{L,1},\cdots,b_{L,m},b^{'}_{L,1},\cdots,b^{'}_{L,m}\}$. So by \cite{BCS},
$\mathbf{\Sigma_{\theta}}=(N_{L},\Sigma_{\theta},\beta_{L})$ is a stacky fan.
\end{pf}

\begin{defn}\label{lawrencetoricdmstack}
The toric DM stack $\mathcal{X}(\mathbf{\Sigma_{\theta}})$ is called the Lawrence toric DM stack.
\end{defn}

For  the map $\beta_{L}^{\vee}:
\mathbb{Z}^{m}\oplus \mathbb{Z}^{m}\to DG(\beta)$ given by 
$(\beta^{\vee},-\beta^{\vee})$, there is an exact sequence
\begin{equation}\Label{exact3}
0\longrightarrow N_{L}^{*}\longrightarrow
\mathbb{Z}^{2m}\stackrel{\beta_{L}^{\vee}}{\longrightarrow}
DG(\beta)\longrightarrow Coker(\beta_{L}^{\vee})\longrightarrow
0.
\end{equation}
Applying $Hom_\mathbb{Z}(-,\mathbb{C}^\times)$ to (\ref{exact3}) yields
\begin{equation}\Label{exact5}
1\longrightarrow \mu\longrightarrow
G\stackrel{\alpha^{L}}{\longrightarrow}
(\mathbb{C}^{\times})^{2m}\longrightarrow
T_{L}\longrightarrow 1, \
\end{equation}
where
$\mu:=Hom_{\mathbb{Z}}(Coker(\beta_{L}^{\vee}),\mathbb{C}^{\times})$ and $T_{L}$ is the torus 
of dimension $m+d$. 
From \cite{BCS} and Proposition \ref{lawrencestacky}, the toric DM stack $\mathcal{X}(\mathbf{\Sigma_{\theta}})$ is
the quotient stack $[(\mathbb{C}^{2m}\setminus V(\mathcal{I}_{\theta}))/G]$, where $G$ acts 
on $\mathbb{C}^{2m}\setminus V(\mathcal{I}_{\theta})$ through the map $\alpha^{L}$ in (\ref{exact5}). 

\subsection*{Hypertoric DM stacks}
Define an ideal in $\mathbb{C}[z,w]$ by:
\begin{equation}\Label{ideal1}
I_{\beta^{\vee}}:=\left\langle\sum_{i=1}^{m}(\beta^{\vee})^{*}(x)_{i}z_{i}w_{i}|~ x\in DG(\beta)^{*}\right\rangle,
\end{equation}
where $(\beta^{\vee})^{*}$ is the map in (\ref{exact1}) and $(\beta^{\vee})^{*}(x)_{i}$ is the $i$-th component
of the vector $(\beta^{\vee})^{*}(x)$.

According to Section 6 in \cite{HS}, $I_{\beta^{\vee}}$ is a prime ideal. Let $Y$ be the closed  subvariety of
$\mathbb{C}^{2m}\setminus V(\mathcal{I}_{\theta})$ determined by the ideal (\ref{ideal1}).
Since 
$(\mathbb{C}^{\times})^{2m}$ acts on $Y$ naturally and the
group $G$ acts on $Y$ through the map $\alpha^{L}$, we have the quotient stack $[Y/G]$. 
Since $Y\subseteq\mathbb{C}^{2m}\setminus V(\mathcal{I}_{\theta})$ is a closed subvariety, the quotient stack $[Y/G]$ is a closed substack 
of $\mathcal{X}(\mathbf{\Sigma_{\theta}})$, and is Deligne-Mumford.

\begin{defn}
The hypertoric Deligne-Mumford stack
$\mathcal{M}(\mathcal{A})$ associated to the  stacky
hyperplane arrangement  $\mathcal{A}$ is defined to be the quotient stack
$[Y/G]$.
\end{defn}

\begin{example}
Let $N=\mathbb{Z}\oplus \mathbb{Z}_{2}$, $\Sigma$ the fan of
projective line $\mathbb{P}^{1}$, and $\beta:
\mathbb{Z}^{3}\to N$  given by $\{b_{1}=(1,0),
b_{2}=(-1,1), b_{3}=(1,0)\}$. Then the Gale dual $\beta^{\vee}:
\mathbb{Z}^{3}\to \mathbb{Z}^{2}$ is given by the
matrix $\left[
\begin{array}{ccc}
1&0&1\\
2&2&0
\end{array}
\right]$. Choose a generic element $\theta=(1,1)$
in $\mathbb{Z}^{2}$ which determines the fan $\Sigma$.
The stacky hyperplane arrangement is $\mathcal{A}=(N,\beta,\theta)$, 
$G=(\mathbb{C}^{\times})^{2}$ and $Y$ is the
subvariety of
$\text{Spec}(\mathbb{C}[z_{1},z_{2},z_{3},w_{1},w_{2},w_{3}])$ determined
by the ideal $I_{\beta^{\vee}}=(z_{1}w_{1}+z_{3}w_{3},
2z_{1}w_{1}+2z_{2}w_{2})$.  Then by \cite{HS}, the coarse moduli space
is the {\em crepant resolution} of the Gorenstein orbifold
$[\mathbb{C}^{2}/\mathbb{Z}_{3}]$, see Figure 3. The corresponding
hyperplane arrangement $\mathcal{H}$ consists of three distinct
points on the real line $\mathbb{R}^{1}$, and the bounded
polyhedron is two segments intersecting at one point. So the core
of the hypertoric variety is two $\mathbb{P}^{1}$ intersecting at
one point.  The hypertoric DM stack $\mathcal{M}(\mathcal{A})$
is a nontrivial $\mu_{2}$-gerbe over the {\em crepant resolution}
according to the action given by the inverse of the above matrix.
If we change $b_{2}$ to $(-1,0)$, we will see an example in
Section \ref{substack} that the hypertoric DM stack is a trivial
$\mu_{2}$-gerbe over the crepant resolution.
\end{example}

\section{Properties of Hypertoric DM Stacks}\label{hypertoricDM}

\subsection*{The coarse moduli space}
Each Deligne-Mumford stack has an underlying coarse moduli space. In this section we prove that the 
coarse moduli space of $\mathcal{M}(\mathcal{A})$ is the underlying 
hypertoric variety.

Consider again the map $\beta^{\vee}: \mathbb{Z}^{m}\rightarrow DG(\beta)$ in (\ref{exact2}),
which is given by the vectors $(a_1,\cdots,a_m)$. As in Section 2, let 
$\overline{\theta}$
be the natural image of $\theta$ under the projection $DG(\beta)\rightarrow \overline{DG(\beta)}$. 
Then the map $\overline{\beta}^{\vee}: \mathbb{Z}^{m}\rightarrow \overline{DG(\beta)}$ is given by $\overline{\beta}^{\vee}=(\overline{a}_{1},\cdots,\overline{a}_{m})$. 
From the map $\overline{\beta}^{\vee}$ we get the  simplicial fan $\Sigma_{\theta}$ in 
(\ref{fan}). By \cite{BCS}, the toric variety $X(\Sigma_{\theta})$, which is  the
geometric quotient $(\mathbb{C}^{2m}-V(\mathcal{I}_{\theta}))\slash G$, is the coarse moduli space of 
the Lawrence toric DM stack $\mathcal{X}(\mathbf{\Sigma_{\theta}})$. The toric variety  $X(\Sigma_{\theta})$
is semi-projective, but not projective. 
In \cite{HS},
from $\beta^{\vee}$ and $\theta$, the authors define the hypertoric variety 
$Y(\beta^{\vee},\theta)$ as the complete intersection of the toric variety 
$X(\Sigma_{\theta})$ by the ideal (\ref{ideal1}), which is the geometric quotient 
$Y\slash G$. We have the following Proposition.

\begin{prop}
The coarse moduli space of $\mathcal{M}(\mathcal{A})$ is
$Y(\beta^{\vee},\theta)$.
\end{prop}
\begin{pf}
Let $X=(\mathbb{C}^{2m}-V(\mathcal{I}_{\theta}))$. By the universal
property of geometric quotients (\cite{KM}), we have that
$X\times_{X(\Sigma_{\theta})}Y(\beta^{\vee},\theta)=Y$. From Lemma
3.3 in \cite{JT}, the stabilizers of points in $X$ are the same as
the stabilizers of the points in $Y$, which are determined by the
box elements in the Lawrence simplicial fan and extended stacky fan.
So we have the following diagram
$$
\vcenter{\xymatrix{
~\mathcal{M}(\mathcal{A})\dto\ar@{^{(}->}[r]\ar@{}[dr]|{\Box}&~\mathcal{X}(\mathbf{\Sigma_{\theta}})\dto \\
Y(\beta^{\vee},\theta) \ar @{^{(}->}[r]& X(\Sigma_{\theta}),}}
$$
which is cartesian. The Lawrence toric variety $X(\Sigma_{\theta})$
is the coarse moduli space of the Lawrence toric Deligne-Mumford
stack $\mathcal{X}(\mathbf{\Sigma_{\theta}})$. So
$\mathcal{M}(\mathcal{A})$ has coarse moduli space
$Y(\beta^{\vee},\theta)$.
\end{pf}
\begin{rmk}
In \cite{HS}, the authors began with the map $\beta^{\vee}$,  and assumed that $DG(\beta)$ and $N$ are free. In our case 
$DG(\beta)$ is a finitely generated abelian group, the toric variety $X(\Sigma_{\theta})$ is again semi-projective
since $\Sigma_{\theta}$ is a semi-projective fan. The hypertoric 
variety $Y(\beta^{\vee},\theta)$ is the complete intersection determined by the ideal (\ref{ideal1}). This reduces to the case in \cite{HS} when $DG(\beta)$ and $N$ are free.
\end{rmk}

\subsection*{Independence of coorientations of hyperplanes}
From (\ref{arrangement}), a hyperplane $H_{i}$ is naturally oriented. Changing the orientation of $H_i$
means changing the map $\beta$ by replacing $b_{i}$ by $-b_{i}$.

\begin{prop}\label{orientation}
$\mathcal{M}(\mathcal{A})$ is independent to the
coorientations of the  hyperplanes in the hyperplane arrangement
$\mathcal{H}=(H_{1},\cdots,H_{m})$ corresponding to the stacky hyperplane arrangement
$\mathcal{A}$.
\end{prop}
\begin{rmk}
Note that changing coorientations {\em does} change the corresponding
normal fan of the weighted polytope $\mathbf{\Gamma}$.
\end{rmk}
\begin{pf}
It  suffices to prove the Proposition when  we change the coorientation of one
hyperplane, say $H_{j}$ for some $j$. Let
$\mathcal{H}'=(H_{1},\cdots,H_{j}',\cdots,H_{m})$. Then we
have a new stacky hyperplane arrangement
$\mathcal{A'}=(N,\beta',\theta)$, where $\beta':
\mathbb{Z}^{m}\to N$ is given by
$\{b_{1},\cdots,-b_{j},\cdots,b_{m}\}$.  Using the technique of
Gale dual in \cite{BCS}, it is easy to check that if the Gale
dual $\beta^{\vee}$ is given by the integral vectors $\beta^{\vee}=(a_{1},\cdots,a_{m})$, then the
Gale dual $(\beta')^{\vee}$ is given by the integral vectors
$(\beta^{'})^{\vee}=(a_{1},\cdots,-a_{j},\cdots,a_{m})$.  Let $\psi:
\mathbb{Z}^{m}\to \mathbb{Z}^{m}$ be the map given by
$e_{i}\mapsto e_{i}$ if $i\neq j$ and $e_{j}\mapsto
-e_{j}$, then we have the following commutative diagrams:
$$\xymatrix{
~\mathbb{Z}^{m}\dto_{\beta}\rto^{\psi} &~\mathbb{Z}^{m}\dto^{\beta^{'}}\\
N\rto^{id} &N,}\quad\quad
\xymatrix{
(\mathbb{Z}^{m})^{*}\dto_{\beta'^{\vee}}\rto^{} &~(\mathbb{Z}^{m})^{*}\dto^{\beta^{\vee}}\\
DG(\beta')\rto^{} &DG(\beta).}
$$
Consider the diagram
$$\xymatrix{
~(\mathbb{Z}^{2m})^{*}\dto_{[\beta'^{\vee},-\beta'^{\vee}]}
\rto^{} &~(\mathbb{Z}^{2m})^{*}\dto^{[\beta^{\vee},-\beta^{\vee}]}\\
DG(\beta')\rto^{} &DG(\beta).}$$ Applying
$Hom_\mathbb{Z}(-,\mathbb{C}^{\times})$ yields the following diagram of
abelian groups
\begin{equation}\Label{diagram1}
\vcenter{\xymatrix{
G\dto_{\alpha^{L}}\rto^{\varphi_{1}} &G'\dto^{(\alpha^{L})'}\\
(\mathbb{C}^{\times})^{2m}\rto^{} &(\mathbb{C}^{\times})^{2m}.}}
\end{equation}

Recall that $Y$ is a subvariety of
$\mathbb{C}^{2m}\setminus V(\mathcal{I}_{\theta})$ defined by the ideal
$I_{\beta^{\vee}}$ in (\ref{ideal1}). When we change the coorientation
of $H_{j}$, the ideals do not change, so $Y'=Y$. By
(\ref{diagram1}), the following diagram is Cartesian:
\begin{equation}
\vcenter{\xymatrix{
Y\times G\dto_{(s,t)}\rto^{\varphi_{0}\times\varphi_{1}} &Y'\times G'\dto^{(s,t)}\\
Y\times Y\rto^{\varphi_{0}\times\varphi_{0}} &Y'\times Y'
,}}
\end{equation}
where $\varphi_{0}$ is determined by the map $\psi$. So the
groupoid $Y\times G\toto Y$ is Morita equivalent to the groupoid
$Y'\times G'\toto Y'$. The stack $[Y/G]$ is isomorphic to
the stack $[Y'/G']$, and $\mathcal{M}(\mathcal{A})\cong
\mathcal{M}(\mathcal{A'})$.
\end{pf}

\begin{rmk}
Let $\mathbf{\Sigma}=(N,\Sigma,\beta)$ be the  extended stacky fan
induced by $\mathcal{A}$.
The toric DM stack $\mathcal{X}(\mathbf{\Sigma})$ is the
quotient stack $[Z/G]$, where $Z=(\mathbb{C}^{n}\setminus
V(J_{\Sigma}))\times (\mathbb{C}^{\times})^{m-n}$ as in
\cite{Jiang}, and $J_{\Sigma}$ is the square-free ideal of the fan
$\Sigma$. So every hypertoric DM stack
$\mathcal{M}(\mathcal{A})$ has an associated  toric DM
stack $\mathcal{X}(\mathbf{\Sigma})$ whose simplicial fan is the
normal fan of the bounded polytope $\mathbf{\Gamma}$ in the
hyperplane arrangement $\mathcal{H}$ determined by the 
stacky hyperplane arrangement $\mathcal{A}$. But by Proposition \ref{orientation}, $\mathcal{M}(\mathcal{A})$ does not determine $\mathcal{X}(\mathbf{\Sigma})$.
\end{rmk}

\begin{example}
Consider Figure 1 again. The corresponding toric variety is
$\mathbb{P}^{2}$. If we change the coorientation of the hyperplane
$2$, then the corresponding normal fan $\Sigma$ of
$\mathbf{\Gamma}$ changes.  The resulting  toric variety is a  Hirzebruch surface. So the associated toric DM stacks are different. But the hypertoric DM stacks are the
same.
\end{example}

\section{Substacks of Hypertoric DM Stacks}\label{substack}

In this section we consider substacks of hypertoric DM
stacks. In particular, we determine the inertia stack of a
hypertoric DM stack.

Let $\mathcal{A}=(N,\beta,\theta)$ be a 
stacky hyperplane arrangement  and 
$\mathbf{\Sigma}=(N,\Sigma,\beta)$  the extended stacky fan induced from $\mathcal{A}$. 
Let $\mathcal{M}(\mathcal{A})$ denote the corresponding hypertoric DM stack.
Consider the map $\beta: \mathbb{Z}^{m}\rightarrow N$ given by $\{b_{1},\cdots,b_{m}\}$.
Let $Cone(\beta)$ be a partially
ordered finite set of cones generated by
$\overline{b}_{1},\cdots,\overline{b}_{m}$. The partial ordering is
defined by requiring that $\sigma\prec\tau$ if $\sigma$ is a face of $\tau$. We have the minimum element $\hat{0}$ which is the cone
consisting of the origin. Let $Cone(\overline{N})$ be the set of
all convex polyhedral cones in the lattice $\overline{N}$. Then we
have a map
$$C: Cone(\beta)\longrightarrow Cone(\overline{N}),$$ such
that for any $\sigma\in Cone(\beta)$, $C(\sigma)$ is the
cone in $\overline{N}$. Then $\Delta_{\mathbf{\beta}}:=(C,Cone(\beta))$ is a simplicial multi-fan in the sense of \cite{HM}. 

\subsection*{Closed substacks}
Recall that in Section \ref{hyperplane} we have the fan $\Sigma_{\theta}$ 
for the  Lawrence toric
variety corresponding to $\pm \beta^{\vee}$.  Let
$\Lambda(\mathcal{B})=\{\overline{b}_{L,1},\cdots,\overline{b}_{L,m},\overline{b}'_{L,1},\cdots,\overline{b}'_{L,m}\}
\subset \overline{N}_{L}$
be the Lawrence lifting of $\mathcal{B}=\{\overline{b}_{1},\cdots,\overline{b}_{m}\}\subset \overline{N}$.
We have the following lemma.

\begin{lem}\label{Lawrence}
If
$\sigma_{\theta}=(\overline{b}_{L,i_{1}},\cdots,\overline{b}_{L,i_{k}},\overline{b}'_{L,i_{1}},
\cdots,\overline{b}'_{L,i_{k}})$ forms a cone in
$\Sigma_{\theta}$, then $\sigma=(\overline{b}_{i_{1}},\cdots,\overline{b}_{i_{k}})$
forms a cone in $\Delta_{\mathbf{\beta}}$.
\end{lem}
\begin{pf}
This can be easily proved from the definition of fan $\Sigma_{\theta}$ in (\ref{fan}).
\end{pf}

For a cone
$\sigma$ in the multi-fan $\Delta_{\mathbf{\beta}}$, let
$link(\sigma)=\{b_{i}: \rho_{i}+\sigma~ \text{is a cone in}~
\Delta_{\mathbf{\beta}}\}$. Then we have a quotient extended
stacky fan
$\mathbf{\Sigma/\sigma}=(N(\sigma),\Sigma/\sigma,\beta(\sigma))$,
where $\beta(\sigma):\mathbb{Z}^{l}\to N(\sigma)$
is given by the images of $\{b_{i}\}$'s in $link(\sigma)$. 
Let $s:=|\sigma|$, then $dim(N_{\sigma})=s$ since $\sigma$ is simplicial.
Consider the commutative diagrams
\[
\begin{CD}
0 @ >>>\mathbb{Z}^{l+s}@ >>> \mathbb{Z}^{m}@ >>>
\mathbb{Z}^{m-l-s} @
>>> 0\\
&& @VV{\widetilde{\beta}}V@VV{\beta}V@VV{}V \\
0@ >>> N @ >{\cong}>>N@ >>> 0 @>>> 0,
\end{CD}
\]
and 
\[
\begin{CD}
0 @ >>>\mathbb{Z}^{s}@ >>> \mathbb{Z}^{l+s}@ >>>
\mathbb{Z}^{l} @
>>> 0\\
&& @VV{\beta_{\sigma}}V@VV{\widetilde{\beta}}V@VV{\beta(\sigma)}V \\
0@ >>> N_{\sigma} @ >{}>>N@ >>> N(\sigma) @>>> 0.
\end{CD}
\]
Applying the
Gale dual yields
\begin{equation}\label{close-hypertoric1}
\begin{CD}
0 @ >>>\mathbb{Z}^{m-l-s}@ >>> \mathbb{Z}^{m}@ >>>
\mathbb{Z}^{l+s} @
>>> 0\\
&& @VV{\cong}V@VV{\beta^{\vee}}V@VV{\widetilde{\beta}^{\vee}}V \\
0@ >>> \mathbb{Z}^{m-l-s} @ >{}>>DG(\beta)@ >{\phi_{1}}>> DG(\widetilde{\beta})@>>>0,
\end{CD}
\end{equation}
and 
\begin{equation}\label{close-hypertoric2}
\begin{CD}
0 @ >>>\mathbb{Z}^{l}@ >>> \mathbb{Z}^{l+s}@ >>>
\mathbb{Z}^{s} @
>>> 0\\
&& @VV{\beta(\sigma)^{\vee}}V@VV{\widetilde{\beta}^{\vee}}V@VV{\beta_{\sigma}^{\vee}}V \\
0@ >>> DG(\beta(\sigma)) @ >{\phi_{2}}>>DG(\widetilde{\beta})@ >>> DG(\beta_{\sigma})
@>>> 0.
\end{CD}
\end{equation}
Since $\mathbb{Z}^{s}\cong N_{\sigma}$, the Gale dual $DG(\beta_{\sigma})=0$. 
And again applying the $Hom_{\mathbb{Z}}(-,\mathbb{C}^{\times})$ functor to the above
two diagrams (\ref{close-hypertoric1}), (\ref{close-hypertoric2}) yields
\begin{equation}\label{close-hypertoric3}
\begin{CD}
1 @ >>>\widetilde{G}@ >>> G@ >>>
(\mathbb{C}^{\times})^{m-l-s} @
>>> 1\\
&& @VV{\widetilde{\alpha}}V@VV{\alpha}V@VV{\cong}V \\
1@ >>> (\mathbb{C}^{\times})^{l+s} @ >{}>>(\mathbb{C}^{\times})^{m}@ >>> (\mathbb{C}^{\times})^{m-l-s}
@>>> 1,
\end{CD}
\end{equation}
and
\begin{equation}\label{close-hypertoric4}
\begin{CD}
1 @ >>>1@ >>> \widetilde{G}@ >{\cong}>>
G(\sigma) @
>>> 1\\
&& @VV{}V@VV{\widetilde{\alpha}}V@VV{\alpha(\sigma)}V \\
1@ >>> (\mathbb{C}^{\times})^{s} @ >{}>>(\mathbb{C}^{\times})^{l+s}@ >>> (\mathbb{C}^{\times})^{l}
@>>> 1.
\end{CD}
\end{equation}
Since $\theta\in DG(\beta)$, from the map $\phi_{1}$ in (\ref{close-hypertoric1}) 
$\theta$ induces $\widetilde{\theta}$ in $DG(\widetilde{\beta})$. From the isomorphism $\phi_{2}$
in (\ref{close-hypertoric2}), we get $\theta(\sigma)\in DG(\beta(\sigma))$. 
Then $\mathcal{A}=(N,\beta,\theta)$ gives 
$\mathcal{A}(\sigma)=(N(\sigma),\beta(\sigma),\theta(\sigma))$ whose induced  
extended stacky fan is $\mathbf{\Sigma\slash \sigma}$.

From (\ref{close-hypertoric1}),(\ref{close-hypertoric2}) we have the following diagrams
\begin{equation}\label{close-hypertoric5}
\vcenter{\xymatrix{
~\mathbb{Z}^{2m}\dto_{(\beta^{\vee},-\beta^{\vee})}\rto^{} &~\mathbb{Z}^{2(l+s)}\dto^{(\widetilde{\beta}^{\vee},-\widetilde{\beta}^{\vee})}\\
DG(\beta)\rto^{} &DG(\widetilde{\beta}) ,}}~~~~~~~~~~~
\vcenter{\xymatrix{
~\mathbb{Z}^{2l}\dto_{(\beta(\sigma)^{\vee},-\beta(\sigma)^{\vee})}\rto^{} &~\mathbb{Z}^{2(l+s)}\dto^{(\widetilde{\beta}^{\vee},-\widetilde{\beta}^{\vee})}\\
DG(\beta(\sigma))\rto^{} &DG(\widetilde{\beta}) .}}
\end{equation}
Taking $Hom_\mathbb{Z}(-,\mathbb{C}^{\times})$ gives
\begin{equation}\label{close-hypertoric6}
\vcenter{\xymatrix{
~\widetilde{G}\dto_{\widetilde{\alpha}^{L}}\rto^{\varphi_{1}} & G\dto^{\alpha^{L}}\\
(\mathbb{C}^{\times})^{2(l+s)} \rto^{}
&(\mathbb{C}^{\times})^{2m} ,}}~~~~~~~~~
\vcenter{\xymatrix{
~\widetilde{G}\dto_{\widetilde{\alpha}^{L}}\rto^{\cong} & G(\sigma)\dto^{\alpha(\sigma)^{L}}\\
(\mathbb{C}^{\times})^{2(l+s)} \rto^{}
&(\mathbb{C}^{\times})^{2l} .}}
\end{equation}
Let $X(\sigma):=(\mathbb{C}^{2l}\setminus V(\mathcal{I}_{\theta(\sigma)}))$ and $Y(\sigma)$ the closed subvariety of 
$X(\sigma)$ defined by the ideal
\begin{equation}\label{idealclose-hypertoric}
I_{\beta(\sigma)^{\vee}}:=\{\sum_{i=1}^{l}(\beta(\sigma)^{\vee})^{*}(x)_{i}z_{i}w_{i}:
\forall x\in DG(\beta(\sigma))^{*}\},
\end{equation}
where $(\beta(\sigma)^{\vee})^{*}: DG(\beta(\sigma))^{*}\rightarrow \mathbb{Z}^{l}$ is 
the dual map of $\beta(\sigma)^{\vee}$ and $(\beta(\sigma)^{\vee})^{*}(x)_{i}$ the $i$-th component of the vector $(\beta(\sigma)^{\vee})^{*}(x)$. Then from the definition of hypertoric DM stacks, we have
$\mathcal{M}(\mathcal{A}(\sigma))=[Y(\sigma)/G(\sigma)]$.
We have the following result:

\begin{prop}\label{closedsub}
If $\sigma$ is a cone in the multi-fan $\Delta_{\mathbf{\beta}}$,
then
$\mathcal{M}(\mathcal{A}(\sigma))$ is a closed substack of
$\mathcal{M}(\mathcal{A})$.
\end{prop}

\begin{pf}
Let $\mathcal{I}_{\theta}$ be the irrelevant ideal in
(\ref{irrelevant}). The hypertoric stack
$\mathcal{M}(\mathcal{A})$ is the quotient stack $[Y/G]$,
where $Y\subset X:=(\mathbb{C}^{2m}\setminus V(\mathcal{I}_{\theta}))$ is the subvariety determined by the
ideal $I_{\beta^{\vee}}$ in (\ref{ideal1}). 

Taking duals to (\ref{close-hypertoric1}), (\ref{close-hypertoric2}) we get:
\begin{equation}\label{close-hypertoric7}
\begin{CD}
0 @ >>>DG(\widetilde{\beta})^{*}@ >>> DG(\beta)^{*}@ >>>
\mathbb{Z}^{m-l-s} @
>>> 0\\
&& @VV{(\widetilde{\beta}^{\vee})^{*}}V@VV{(\beta^{\vee})^{*}}V@VV{\cong}V \\
0@ >>> \mathbb{Z}^{l+s} @ >{}>>\mathbb{Z}^{m}@ >{\phi_{1}}>> \mathbb{Z}^{m-l-s}@>>>0,
\end{CD}
\end{equation}
and 
\begin{equation}\label{close-hypertoric8}
\begin{CD}
0 @ >>>0@>>>DG(\widetilde{\beta})^{*}@ >{\cong}>> DG(\beta(\sigma))^{*}@ 
>>> 0\\
&& @VV{}V@VV{(\widetilde{\beta}^{\vee})^{*}}V@VV{(\beta(\sigma)^{\vee})^{*}}V \\
0@ >>> \mathbb{Z}^{s} @ >{}>>\mathbb{Z}^{l+s}@ >{\phi_{1}}>> \mathbb{Z}^{l}@>>>0,
\end{CD}
\end{equation}
Let $W(\sigma)$ be the subvariety of $X$ defined
by the ideal $J(\sigma):=\left<z_{i},w_{i}:\rho_{i}\subseteq
\sigma\right>$. Then $W(\sigma)$ contains the $\mathbb{C}$-points
$(z,w)\in \mathbb{C}^{2m}$ such that 
the cone spanned by $\{\rho_{i}: z_{i}=w_{i}=0\}$ containing 
$\sigma$ belongs to $\Delta_{\mathbf{\beta}}$. From Lemma \ref{Lawrence},
the $\mathbb{C}$-point $(z,w)$ in $W(\sigma)$ such that 
$\rho_{i}\nsubseteq \sigma\cup link(\sigma)$ implies that 
$z_{i}\neq 0$ or $w_{i}\neq 0$. It is clear that $W(\sigma)$ is invariant under the $G$-action defined by (\ref{exact5}). 
Let $V(\sigma):=Y\cap W(\sigma)$. Then from  (\ref{ideal1}),(\ref{idealclose-hypertoric}) and 
(\ref{close-hypertoric7}),(\ref{close-hypertoric8}),
$V(\sigma)\cong Y(\sigma)\times(\mathbb{C}^{\times})^{m-s-l}\times (0)^{m-s-l}$ and the components
$0$ are determined by the choice of the generic element $\theta$.

Let $\varphi_{0}: Y(\sigma)\longrightarrow V(\sigma)$ be the inclusion 
given by $(z,w)\longmapsto (z,w,1,0)$. From the map $\varphi_{1}$ in (\ref{close-hypertoric6}), we have 
a morphism of groupoids $\varphi_{0}\times \varphi_{1}: Y(\sigma)\times G(\sigma)\toto V(\sigma)\times G$ which 
induces a morphism of stacks $\varphi: [Y(\sigma)/G(\sigma)]\longrightarrow [V(\sigma)/G]$.
To prove that it is an isomorphism, we first prove that the following 
diagram is cartesian:
$$\xymatrix{
Y(\sigma)\times G(\sigma)\dto_{(s,t)}\rto^{\varphi_{0}\times\varphi_{1}} &V(\sigma)\times G\dto^{(s,t)}\\
Y(\sigma)\times Y(\sigma)\rto^{\varphi_{0}\times\varphi_{0}}
&V(\sigma)\times V(\sigma).} $$ 
This is easy to prove. Given an element $((z_1,w_1),(z_2,w_2))\in Y(\sigma)\times Y(\sigma)$, under 
the map $\varphi_{0}\times\varphi_{0}$, we get $((z_1,w_1,1,0),(z_2,w_2,1,0))\in V(\sigma)\times V(\sigma)$.
If there is an element $g\in G$ such that $g(z_1,w_1,1,0)=(z_2,w_2,1,0)$, then from the exact sequence in the first row
of (\ref{close-hypertoric3}), there is an element $g(\sigma)\in G(\sigma)$ such that $g(\sigma)(z_1,w_1)=(z_2,w_2)$.
Thus we have an element $((z_1,w_1),g(\sigma))\in Y(\sigma)\times G(\sigma)$. So the 
morphism $\varphi: [Y(\sigma)/G(\sigma)]\longrightarrow [V(\sigma)/G]$ is injective. 
Let $(z,w,s,0)$ be an element in $V(\sigma)$, then there exists an element $g\in (\mathbb{C}^{\times})^{m-l-s}$ such that 
$g(z,w,s,0)=(z,w,1,0)$. From (\ref{close-hypertoric3}), $g$ determines an element in $G$, so
$\varphi$ is surjective and $\varphi$ is an isomorphism.
Clearly  the stack $[V(\sigma)/G]$ is a closed substack of $\mathcal{M}(\mathcal{A})$, so the stack
$\mathcal{M}(\mathcal{A}(\sigma))=[Y(\sigma)/G(\sigma)]$ is
also a closed substack of $\mathcal{M}(\mathcal{A})$.
\end{pf}

\subsection*{Open substacks}
We now study open substacks of $\mathcal{M}(\mathcal{A})$.
Let $\sigma$ be a top dimensional cone in
$\Delta_{\mathbf{\beta}}$. Then $\mathbf{\sigma}
=(\mathbb{Z}^{d},\sigma,\beta_{\sigma})$ is a stacky fan, where
$\beta_{\sigma}: \mathbb{Z}^{d}\to N$ is given by
$b_{i}$ for $\rho_{i}\subseteq \sigma$. 
Since $N$ has rank $d$, we find that $DG(\beta_{\sigma})$ is a finite abelian
group. So in this case the generic element $\theta$ induces zero in  $DG(\beta_{\sigma})$. This is
the degenerate case, which means that the corresponding ideal (\ref{ideal1}) is zero.  Thus 
$$Y_{\sigma}=\mathbb{C}^{2d}.$$
Note that $G_{\sigma}=Hom_{\mathbb{Z}}(DG(\beta_{\sigma}),\mathbb{C}^{\times})$ is a
finite abelian group.
According to the construction of hypertoric
DM stack in Section 2, the hypertoric DM stack
$\mathcal{M}(\mathcal{\sigma})$ associated to $\mathbf{\sigma}$ is the quotient stack
$[Y_{\sigma}/G_{\sigma}]$ 
which can be regarded as a local chart of the hypertoric orbifold $[Y/G]$.

\begin{prop}
If $\sigma$ is a top-dimensional cone in the multi-fan
$\Delta_{\mathbf{\beta}}$, then $\mathcal{M}(\sigma)$ is
an open substack of $\mathcal{M}(\mathcal{A})$.
\end{prop}

\begin{pf}
Since $\sigma$ is a top dimensional cone in $\Delta_{\mathbf{\beta}}$, from
(\ref{exact2}) we get a basis $C$ of $\overline{DG(\beta)}$. Let $U_{\sigma}$ be the open subvariety of
$\mathbb{C}^{2m}\backslash V(\mathcal{I}_{\theta})$ defined by the
monomials $\prod C(\theta)$ in (\ref{irrelevant}). Let $V_{\sigma}=U_{\sigma}\cap
Y$, i.e. the points in $U_{\sigma}$ staisfying (\ref{ideal1}). 
Then we have the groupoid
$V_{\sigma}\times G\toto V_{\sigma}$ associated to the action of
$G$ on $V_{\sigma}$. It is clear that this groupoid defines an
open substack of $\mathcal{M}(\mathcal{A})$. Next we show that
this substack is isomorphic to $\mathcal{M}(\sigma)$.

Consider the following commutative diagram:
\[
\begin{CD}
0 @ >>>\mathbb{Z}^{d}@ >>> \mathbb{Z}^{m}@ >>>
\mathbb{Z}^{m-d} @
>>> 0\\
&& @VV{\beta_{\sigma}}V@VV{\beta}V@VV{\beta^{'}}V \\
0@ >>>N @ >{id}>>N@ >>> 0 @>>> 0.
\end{CD}
\]
Applying Gale dual and $Hom_\mathbb{Z}(-,\mathbb{C}^{\times})$, we obtain 
\begin{equation}\Label{propsecond}
\begin{CD}
1 @ >>>G_{\sigma}@ >{\varphi_{1}}>> G@ >>>
(\mathbb{C}^{\times})^{m-d} @
>>> 1\\
&& @VV{\alpha_{\sigma}}V@VV{\alpha}V@VV{id}V \\
1@ >>>(\mathbb{C}^{\times})^{d} @
>{}>>(\mathbb{C}^{\times})^{m}@ >>> (\mathbb{C}^{\times})^{m-d} @>>> 1.
\end{CD}
\end{equation}
We construct a morphism $\varphi_{0}: Y_{\sigma}\to V_{\sigma}$.
For $\rho_{j}\nsubseteq \sigma$, we set $z_{j}=1$ if $z_{j}$ is component 
of a monomial of $C(\theta)$ in (\ref{irrelevant}) or $w_{j}=1$ if $w_{j}$ is component 
of a monomial of $C(\theta)$ in (\ref{irrelevant}), then from (\ref{ideal1}), 
the corresponding $w_{j}$ or $z_{j}$ can be represented as linear component of 
$\{z_{i}w_{i}\}$ for $\rho_{i}\subseteq\sigma$. 
Let $\widetilde{\varphi}_{0}: \mathbb{C}^{2d}\rightarrow U_{\sigma}$ be the morphism
given by $z_{i},w_{i}\longmapsto z_{i},w_{i}$ for $\rho_{i}\subseteq\sigma$,
and $z_{j},w_{j}$ to the corresponding $1$ or linear combination of $\{z_{i}w_{i}\}$ for 
$\rho_{i}\subseteq\sigma$ in the above analysis. Then $\widetilde{\varphi}_{0}$
induces a morphism $\varphi_{0}: Y_{\sigma}\to V_{\sigma}$.
Hence we have a morphism of groupoids
$$\Phi:=(\varphi_{0}\times \varphi_{0},\varphi_{0}\times \varphi_{1}):
[Y_{\sigma}\times G_{\sigma}\toto Y_{\sigma}]\longrightarrow
[V_{\sigma}\times G\toto V_{\sigma}],$$
where $\varphi_{1}$ is the morphism in (\ref{propsecond}). This morphism determines a
morphism of the associated stacks. The isomorphism of these two
stacks comes from the following Cartesian diagram:
\begin{equation}\Label{prop4.6}
\vcenter{\xymatrix{
Y_{\sigma}\times G_{\sigma}\dto_{(s,t)}\rto^{\varphi_{0}\times\varphi_{1}} &V_{\sigma}\times G\dto^{(s,t)}\\
Y_{\sigma}\times Y_{\sigma}\rto^{\varphi_{0}\times\varphi_{0}}
&V_{\sigma}\times V_{\sigma}.}} 
\end{equation}
\end{pf}
\subsection*{Inertia stacks}
Let $N_{\sigma}$ be the sublattice generated by $\sigma$, and $N(\sigma):=N/N_{\sigma}$. Note that when $\sigma$ is a  top
dimensional cone, $N(\sigma)$ is the local orbifold group in the
local chart of the coarse moduli space of the hypertoric toric DM
stack. Namely:

\begin{lem}\label{finitegroup}
Let $\sigma$ be a top-dimensional cone in the multi-fan
$\Delta_{\mathbf{\beta}}$. Then $G_{\sigma}\cong N(\sigma)$.
\end{lem}

\begin{pf}
The proof is the same as the proof for a top dimensional cone in a
simplicial fan in Proposition 4.3 in \cite{BCS}.
\end{pf}

Recall that $G$ acts on $(\mathbb{C}^\times)^{2m}$ via the map $\alpha^{L}:
G\to (\mathbb{C}^{\times})^{2m}$ in (\ref{exact5}). We write 
$$\alpha^{L}(g)=(\alpha^{L}_{1}(g),\cdots,\alpha^{L}_{m}(g),\alpha^{L}_{1+m}(g),\cdots,\alpha^{L}_{2m}(g)).$$ 

\begin{lem}\label{fixpoint}
Let $(z,w)\in Y$ be a point fixed by $g\in G$. If $\alpha_i^L(g)\neq 1$, then $z_{i}=w_{i}=0$.
\end{lem}

\begin{pf}  
Since $G$ acts on $\mathbb{C}^{2m}$ through the matrix $\beta_{L}^{\vee}=[\beta^{\vee},-\beta^{\vee}]$ in
(\ref{exact3}), we have that
$\alpha^{L}_{i+m}(g)=\alpha^{L}_{i}(g)^{-1}$. The Lemma follows immediately.
\end{pf}

Given the multi-fan $\Delta_{\mathbf{\beta}}$, we consider the pairs
$(v,\sigma)$, where $\sigma$ is a cone in $\Delta_{\mathbf{\beta}}$, 
$v\in N$ such that $\overline{v}=\sum_{\rho_{i}\subseteq \sigma}\alpha_{i}b_{i}$ for 
$0<\alpha_{i}<1$. Note that $\sigma$ is the minimal cone in $\Delta_{\mathbf{\beta}}$
satisfying the above condition. Let $Box(\Delta_{\mathbf{\beta}})$ be the set of all such pairs  $(v,\sigma)$.

\begin{prop}\label{fixbox}
There is an one-to-one correspondence between $g\in G$ with
nonempty fixed point set and $(v,\sigma)\in Box(\Delta_{\mathbf{\beta}})$.
Moreover, for such $g$ and $(v,\sigma)$ we have $[Y^{g}/G]\cong
\mathcal{M}(\mathcal{A}(\sigma))$.
\end{prop}

\begin{pf}
Let $(v,\sigma)\in Box(\Delta_{\mathbf{\beta}})$. Since $\sigma$ is 
contained in a top dimensional cone $\tau$ in $\Delta_{\mathbf{\beta}}$, we have
$v\in N(\tau)$. By Lemma \ref{finitegroup},
$N(\tau)\cong G_{\tau}$. Hence $v$ determines an
element in $G_{\tau}$.  Using the morphism $\varphi_{1}$ in
(\ref{propsecond}), we see that $g$ fixes a point in $Y$.

Conversely, suppose $g\in G$ fixes a point $(z,w)$ in $Y$, where
$(z,w)\in \mathbb{C}^{2m}$. By
Lemma \ref{fixpoint}, the point $(z,w)$ satisfies the condition
that if $\alpha_i^L(g)\neq 1$ then $z_{i}=w_{i}=0$.  From the definition of
$\mathbb{C}^{2m}\backslash V(\mathcal{I}_{\theta})$, there is a
cone in $\Sigma_{\theta}$ containing the rays for which $z_{i}=w_{i}=0$. By Lemma \ref{Lawrence}, the rays $\rho_{i}$ for which
$z_{i}=0$ is a cone in $\Delta_{\mathbf{\beta}}$ which we call $\sigma$. So $g$ stabilizes 
$Y_{\tau}=\mathbb{C}^{2d}$ in $V_{\tau}$ through 
$\varphi_{0}$ in (\ref{prop4.6}) for any top dimensional cone
$\tau$ containing $\sigma$, and $g$ corresponds to an element
$(v,\sigma)\in Box(\Delta_{\mathbf{\beta}})$. From the definition of
$W(\sigma)$ and $V(\sigma)$ in Proposition \ref{closedsub}, we
have $W(\sigma)\cong Y^{g}$ and $[V(\sigma)/G]\cong [Y^{g}/G]$ which is
$\mathcal{M}(\mathcal{A}(\sigma))$.
\end{pf}

We determine the inertia stack of a hypertoric DM stack.

\begin{prop}\label{inertia}
The inertia stack of $\mathcal{M}(\mathcal{A})$ is given by
$$I(\mathcal{M}(\mathcal{A}))=\coprod_{(v,\sigma)\in
Box(\Delta_{\mathbf{\beta}})}\mathcal{M}(\mathcal{A}(\sigma)).$$
\end{prop}
\begin{pf}
The hypertoric DM stack $\mathcal{M}(\mathcal{A})=[Y/G]$ is a
quotient stack. Its  inertia stack is determined as
$$I\left(\mathcal{M}(\mathcal{A})\right)=\left[\left(\coprod_{g_\in G}
Y^{g}\right)/ G\right].$$ By Proposition \ref{fixbox}, the stack
$[Y^{g}/G]$  is isomorphic to the stack
$\mathcal{M}(\mathcal{A}(\sigma))$ for some $(v,\sigma)\in
Box(\Delta_{\mathbf{\beta}})$.
\end{pf}

\begin{example}
Let $\mathbf{\Sigma}=(N,\Sigma,\beta)$ be an extended stacky fan,
where $N=\mathbb{Z}^{2}$, the simplicial
\begin{center}
$\begin{picture}(0,0)%
\includegraphics{figure2.pstex}%
\end{picture}%
\setlength{\unitlength}{3947sp}%
\begingroup\makeatletter\ifx\SetFigFont\undefined%
\gdef\SetFigFont#1#2#3#4#5{%
  \reset@font\fontsize{#1}{#2pt}%
  \fontfamily{#3}\fontseries{#4}\fontshape{#5}%
  \selectfont}%
\fi\endgroup%
\begin{picture}(3858,2364)(3589,-4465)
\end{picture}
$
\end{center}
fan $\Sigma$ is the fan of weighted projective plane
$\mathbb{P}(1,2,2)$, and $\beta: \mathbb{Z}^{4}\to N$
is given by the vectors
$\{b_{1}=(1,0),b_{2}=(0,1),b_{3}=(-2,-2),b_{4}=(0,-1)\}$, where
$b_{1},b_{2},b_{3}$ are the generators of the rays in $\Sigma$.  
Choose generic element $\theta=(1,1)\in DG(\beta)\cong
\mathbb{Z}^{2}$. Then $\mathcal{A}=(N,\beta,\theta)$ is the stacky hyperplane arrangement 
whose induced extended stacky fan is $\mathbf{\Sigma}$. 
A lifting of $\theta$ in $\mathbb{Z}^{4}$ through
the Gale dual map $\beta^{\vee}$ is $r=(1,1,-3,0)$. The
corresponding hyperplane arrangement
$\mathcal{H}=(H_{1},H_{2},H_{3},H_{4})$ consists of 4 lines, see
Figure 2. 
Take $v=\frac{1}{2}b_{3}$, then $(v,\sigma)\in
Box(\Delta_{\mathbf{\beta}})$, where  $\sigma$ is the ray generated
by $b_{3}$. Consider the following diagram
\[
\begin{CD}
0 @ >>>\mathbb{Z}@ >>> \mathbb{Z}^{4}@ >>> \mathbb{Z}^{3} @
>>> 0\\
&& @VV{\beta_{\sigma}}V@VV{\beta}V@VV{\beta(\sigma)}V \\
0@ >>> N_{\sigma} @ >{}>>N@ >>> \mathbb{Z}\oplus\mathbb{Z}_{2}
@>>> 0.
\end{CD}
\]
We have the quotient extended stacky fan
$\mathbf{\Sigma/\sigma}=(N(\sigma),\Sigma/\sigma,\beta(\sigma))$,
where $\beta(\sigma): \mathbb{Z}^{3}\to N(\sigma)$ is
given by the vectors $\{(1,0),(-1,0),(1,0)\}$, and $(1,0)$ is the
extra data in the quotient extended stacky fan. Taking Gale dual, we get
\[
\begin{CD}
0 @ >>>\mathbb{Z}^{3}@ >>> \mathbb{Z}^{4}@ >>> \mathbb{Z} @
>>> 0\\
&& @VV{\beta(\sigma)^{\vee}}V@VV{\beta^{\vee}}V@VV{\beta_{\sigma}^{\vee}}V \\
0@ >>> \mathbb{Z}^{2}\oplus\mathbb{Z}_{2} @ >{}>>\mathbb{Z}^{2}@
>>> 0@>>> 0,
\end{CD}
\]
where $\beta^{\vee}$ is given by the matrix $\left[
\begin{array}{cccc}
2&2&1&0\\
0&1&0&1
\end{array}
\right]$ and $\beta(\sigma)^{\vee}$ is given by $\left[
\begin{array}{ccc}
2&2&1\\
0&1&0
\end{array}
\right]$. The associated generic element $\theta(\sigma)=(1,1,0)$
and the lifting of $\theta(\sigma)$ in $\mathbb{Z}^{3}$ is
$r(\sigma)=(1,1,-3)$. So the quotient hyperplane arrangement  $\mathcal{A}(\sigma)=(N(\sigma),\beta(\sigma),\theta(\sigma))$
is a line with three distinct points
$\{-1,1,3\}$. The bounded polyhedron of this hyperplane
arrangement is two segments intersecting at one point, see Figure
3.
\begin{center}
\begin{picture}(0,0)%
\includegraphics{xfig_tutoring.pstex}%
\end{picture}%
\setlength{\unitlength}{3947sp}%
\begingroup\makeatletter\ifx\SetFigFont\undefined%
\gdef\SetFigFont#1#2#3#4#5{%
  \reset@font\fontsize{#1}{#2pt}%
  \fontfamily{#3}\fontseries{#4}\fontshape{#5}%
  \selectfont}%
\fi\endgroup%
\begin{picture}(2615,591)(3751,-3865)
\end{picture}

\end{center}
The core of $\mathcal{M}(\mathcal{A}(\sigma))$ corresponds to these two
segments, hence  is two $\mathbb{P}^{1}$'s meeting at one
point. Adding the stacky structure  the twisted sector
$\mathcal{M}(\mathcal{A}(\sigma))$ corresponding to the element
$v$ is the trivial $\mu_{2}$-gerbe over the  {\em crepant resolution}
of the stack $[\mathbb{C}^{2}/\mathbb{Z}_{3}]$.
\end{example}

\section{Orbifold Chow ring of $\mathcal{M}(\mathcal{A})$}\label{chow}

In this section we discuss the orbifold Chow ring of hypertoric DM
stacks. We determine its module structure, then compute the orbifold cup product. 

\subsection{The module structure}
We first consider the ordinary Chow ring for hypertoric DM stacks.
According to \cite{K}, the cohomology ring of
$\mathcal{M}(\mathcal{A})$ is generated by the Chern classes
of some line bundles defined as follows. Applying
$Hom_\mathbb{Z}(-,\mathbb{C}^\times)$ to (\ref{exact2}), we have
$$1\longrightarrow \mu\longrightarrow
G\stackrel{\alpha}{\longrightarrow}
(\mathbb{C}^{\times})^{m}\longrightarrow T\longrightarrow 1.$$

\begin{defn}\label{linebdle}
For every $b_{i}$ in the stacky hyperplane arrangement, define the line
bundle $L_{i}$ over $\mathcal{M}(\mathcal{A})$ to be the trivial line bundle $Y\times \mathbb{C}$ with the $G$-action on $\mathbb{C}$ defined via the $i$-th component of the morphism $\alpha: G\to
(\mathbb{C}^{\times})^{m}$ in the above exact sequence.
\end{defn}

For any $c\in N$, there is a cone 
$\sigma\in \Delta_\mathbf{\beta}$ such that 
$\overline{c}=\sum_{\rho_{i}\subseteq \sigma}\alpha_{i}\overline{b}_{i}$ where 
$\alpha_{i}>0$ are  rational numbers. Let $N^{\Delta_\mathbf{\beta}}$ denote  
all the pairs $(c,\sigma)$. Then $N^{\Delta_\mathbf{\beta}}$ gives rise a 
group ring
$$\mathbb{Q}[\Delta_\mathbf{\beta}]=\bigoplus_{(c,\sigma)\in N^{\Delta_\mathbf{\beta}}}\mathbb{Q}\cdot y^{(c,\sigma)},$$
where $y$ is a formal variable. By abuse of notation, we write $y^{(b_{i},\rho_i)}$ as $y^{b_{i}}$.
The multiplication is given in terms of the ceiling function for fans which we define below.
Since the multi-fan $\Delta_\mathbf{\beta}$ is simplicial, we have the following Lemma.
\begin{lem}\Label{lemma4.5}
For any $c\in N$, there exists a unique cone  $\sigma\in\Delta_{\mathbf{\beta}}$
and $(v,\tau)\in Box(\Delta_{\mathbf{\beta}})$ such that $\tau\subseteq\sigma$ and 
$$c=v+\sum_{\rho_{i}\subseteq\sigma}m_{i}b_{i}$$
where $m_{i}\in\mathbb{Z}_{\geq 0}$. ~$\square$
\end{lem}

\begin{defn}\Label{fraction}
$(v,\tau)$ is called the fractional part of $(c,\sigma)$.
\end{defn}

Now for $(c,\sigma)\in N^{\Delta_\mathbf{\beta}}$, from Lemma \ref{lemma4.5}, we write 
$c=v+\sum_{\rho_{i}\subseteq \sigma}m_{i}b_{i}$, where $m_{i}$'s are nonnegative integers. 
We define  the {\em ceiling function} $\lceil c \rceil_{\sigma}$ by  
$$\lceil c \rceil_{\sigma}=\sum_{\rho_{i}\subseteq \tau}b_{i}+\sum_{\rho_{i}\subseteq \sigma}m_{i}b_{i}.$$
Note that if $\overline{v}=0$, $\lceil c \rceil_{\sigma}=\sum_{\rho_{i}\subseteq \sigma}m_{i}b_{i}$.  
For two pairs  $(c_1,\sigma_1)$, $(c_2,\sigma_2)$, if $\sigma_{1}\cup\sigma_{2}$ is a cone in $\Delta_\mathbf{\beta}$, define 
$\epsilon(c_1,c_2):=\lceil c_1 \rceil_{\sigma_{1}}+\lceil c_2 \rceil_{\sigma_{2}}-\lceil c_1+c_2 \rceil_{\sigma_{1}\cup\sigma_2}$.
Let $\sigma_{\epsilon}\subseteq\sigma_1\cup\sigma_2$ be the minimal cone in $\Delta_\mathbf{\beta}$ containing $\epsilon(c_1,c_2)$ so that 
$(\epsilon(c_1,c_2),\sigma_{\epsilon})\in N^{\Delta_\mathbf{\beta}}$.  The ceiling function $\lceil c \rceil_{\sigma}$ is an integral linear combination 
of $b_{i}$'s for $\rho_{i}\subseteq\sigma$. We define the grading on $\mathbb{Q}[\Delta_\mathbf{\beta}]$ as follows.
For any  $(c,\sigma)$, write $c=v+\sum_{\rho_{i}\subseteq \sigma}m_{i}b_{i}$,  then
$$deg(y^{(c,\sigma)}):=|\tau|+\sum_{\rho_{i}\subseteq\sigma}m_{i},$$ where $|\tau|$ is the dimension of $\tau$. Let
$Cir(\Delta_\mathbf{\beta})$ be the ideal in
$\mathbb{Q}[\Delta_\mathbf{\beta}]$ generated by the elements in (\ref{ideal2}). 
The multiplication $y^{(c_{1},\sigma_{1})}\cdot y^{(c_{2},\sigma_{2})}$
is defined by (\ref{product}). 

\begin{lem}\Label{associative}
The multiplication (\ref{product}) is associative.
\end{lem}

\begin{pf}
For any three pairs $(c_1,\sigma_1),(c_2,\sigma_2),(c_3,\sigma_3)$, if 
$\sigma_{1}\cup\sigma_{2}\cup\sigma_3$ is a cone in $\Delta_{\beta}$, let 
$\sigma\subseteq\sigma_1\cup\sigma_2\cup\sigma_3$ be the minimal cone in $\Delta_{\beta}$ containing 
$$\epsilon(c_1,c_2,c_3):=\lceil c_1 \rceil_{\sigma_{1}}+\lceil c_2 \rceil_{\sigma_{2}}+\lceil c_3 \rceil_{\sigma_{3}}-\lceil c_1+c_2+c_3 \rceil_{\sigma_{1}\cup\sigma_2\cup\sigma_3},$$ 
such that $(\epsilon(c_1,c_2,c_3),\sigma)\in N^{\Delta_{\beta}}$. Then we check from the properties of ceiling function that 
$(y^{(c_1,\sigma_1)}\cdot y^{(c_2,\sigma_2)})\cdot y^{(c_3,\sigma_3)}$ and $y^{(c_1,\sigma_1)}\cdot (y^{(c_2,\sigma_2)}\cdot y^{(c_3,\sigma_3)})$ are both equal to 
$$
\begin{cases}
(-1)^{|\sigma|}y^{(c_{1}+c_{2}+c_3+\epsilon(c_1,c_2,c_3),\sigma_{1}\cup\sigma_{2}\cup\sigma_3)}&\text{if
$\sigma_{1}\cup\sigma_{2}\cup\sigma_3$ is a cone in $\Delta_{\mathbf{\beta}}$}\,,\\
0&\text{otherwise}\,.
\end{cases} 
$$
It is easy to check that the product preserves the grading, and the proof is left to the readers.
\end{pf}

Consider the map $\beta: \mathbb{Z}^{m}\rightarrow N$ which is given by $\{b_{1},\cdots,b_{m}\}$. 
We take $\{1,\cdots,m\}$ as the 
vertex set. The {\em matroid complex} $M_{\beta}$ is defined using
$\beta$ by requiring that $F\in M_{\beta}$ iff the normal vectors
$\{\overline{b}_{i}\}_{i\in F}$ are  linearly independent in $\overline{N}$. The
{\em Stanley-Reisner ring}  of the matroid $M_{\beta}$ is
$$\mathbb{Q}[M_{\beta}]=\frac{\mathbb{Q}[y^{b_{1}},\cdots,y^{b_{m}}]}{I_{M_{\beta}}},$$
where $I_{M_{\beta}}$ is the matroid ideal generated by the set
of  square-free monomials
$$\{y^{b_{i_{1}}}\cdots y^{b_{i_{k}}}|
\overline{b}_{i_{1}},\cdots,\overline{b}_{i_{k}} ~\text{linearly
dependent in}~\overline{N}\}.$$  
It is clear that $\mathbb{Q}[M_{\beta}]$ is a subring of 
$\mathbb{Q}[\Delta_\mathbf{\beta}]$ under the injection
$y^{b_{i}}\longmapsto y^{(b_{i},\rho_{i})}$. 
 
\begin{lem}\label{coarsecoh}
Let $\mathcal{A}=(N,\beta,\theta)$ be a stacky hyperplane arrangement
and $\mathcal{M}(\mathcal{A})$ the corresponding hypertoric
DM stack, then we have an isomorphism of graded rings
$$A^{*}(\mathcal{M}(\mathcal{A}))\cong \frac{\mathbb{Q}[M_{\beta}]}{Cir(\Delta_\mathbf{\beta})},$$
given by $c_{1}(L_{i})\mapsto y^{b_{i}}$, where $Cir(\Delta_\mathbf{\beta})$ is the ideal generated by elements in (\ref{ideal2}).
\end{lem}
\begin{pf}
Let $Y(\beta^{\vee},\theta)$ be the coarse moduli space of the hypertoric DM
stack $\mathcal{M}(\mathcal{A})$. By \cite{HS}, we have
$$A^{*}(Y(\beta^{\vee},\theta))\cong \frac{\mathbb{Q}[M_{\beta}]}{Cir(\Delta_\mathbf{\beta})},$$
 given by $D_{i}\mapsto y^{b_{i}}$, where $D_{i}$ is the
$T$-equivariant Weil divisor on $Y(\beta^{\vee},\theta)$.  Let $a_{i}$ be the
first lattice vector in the ray generated by $b_{i}$, then
$\overline{b}_{i}=l_{i}a_{i}$ for some positive integer $l_{i}$.
By \cite{V},  the Chow ring of the stack
$\mathcal{M}(\mathcal{A})$ is isomorphic to the Chow ring of
its coarse moduli space $Y(A,\theta)$ via $c_{1}(L_{i})\mapsto
l_{i}^{-1}\cdot D_{i}$, and
$\sum_{i=1}^{m}e(a_{i})l_{i}y^{b_{i}}=\sum_{i=1}^{m}e(b_{i})y^{b_{i}}$
for $e\in N^{*}$.
\end{pf}
\begin{rmk}
Note that we have an isomorphism:
$$\frac{\mathbb{Q}[M_{\beta}]}{Cir(\Delta_\mathbf{\beta})}\cong
\frac{\mathbb{Q}[y^{b_{1}},\cdots,y^{b_{m}},y^{b_{1}^{'}},\cdots,y^{b_{m}^{'}}]}
{(Cir(\Delta_\mathbf{\beta}),y^{b_{1}}+y^{b_{1}^{'}},\cdots,y^{b_{m}}+y^{b_{m}^{'}})}$$
for which the right side is the Chow ring of the Lawrence toric DM stack 
$\mathcal{X}(\mathbf{\Sigma_{\theta}})$ defined in Definition \ref{lawrencetoricdmstack}, see \cite{HS}.
So in the DM stack situation $A^{*}(\mathcal{M}(\mathcal{A}))\cong A^{*}(\mathcal{X}(\mathbf{\Sigma_{\theta}}))$.
This point of view will be covered and generalized in \cite{JT}.
\end{rmk}

Let $A_{orb}^*(\mathcal{M}(\mathcal{A}))$ denote the orbifold
Chow ring of $\mathcal{M}(\mathcal{A})$, which by definition
is $A^*(I(\mathcal{M}(\mathcal{A})))$ as a group. By
Proposition \ref{inertia}, we have
$$A^*(I(\mathcal{M}(\mathcal{A})))\cong \bigoplus_{(v,\sigma)\in
Box(\Delta_{\mathbf{\beta}})}A^{*}(\mathcal{M}(\mathcal{A}(\sigma))).$$
For $(v,\sigma)\in Box(\Delta_{\mathbf{\beta}})$, there is an exact sequence of vector bundles,
$$0\to T\mathcal{M}(\mathcal{A}(\sigma))\to
T\mathcal{M}(\mathcal{A})|_{\mathcal{M}(\mathcal{A}(\sigma))}\to N_v\to 0,$$
where $N_v$ denotes the normal bundle of $\mathcal{M}(\mathcal{A}(\sigma))$
in $\mathcal{M}(\mathcal{A})$. On the other hand, there is a surjective morphism
$$\bigoplus_{i=1}^m (L_i\oplus L_i^{-1})\to T\mathcal{M}(\mathcal{A}).$$
Restricting this to $\mathcal{M}(\mathcal{A}(\sigma))$ yields a surjective map
$$\bigoplus_{\rho_i\subset \sigma(\overline{v})} (L_i\oplus L_i^{-1})\to N_v.$$
Moreover, the element in the local group represented by $v$ acts trivially on the kernel.
Let $v$ act on $L_i$ with eigenvalue $e^{2\pi\sqrt{-1}w_i}$, where $w_i\in [0,1)\cap\mathbb{Q}$.
It follows that the age function on the component $\mathcal{M}(\mathcal{A}(\sigma))$
assumes the value $$\sum_{\rho_i\subset \sigma} (w_i+ (1-w_i))=|\sigma|.$$
Hence $A_{orb}^*(\mathcal{M}(\mathcal{A}))$ as a graded module coincides with
$$\bigoplus_{(v,\sigma)\in
Box(\Delta_{\mathbf{\beta}})}A^{*}(\mathcal{M}(\mathcal{A}(\sigma)))[|\sigma|].$$
Note that $A_{orb}^*(\mathcal{M}(\mathcal{A}))$ is
$\mathbb{Z}$-graded, due to the fact that
$\mathcal{M}(\mathcal{A})$ is hyperk\"ahler.

Again since the multi-fan $\Delta_{\mathbf{\beta}}$  is simplicial, we have the following result, similar to 
Lemma 4.6 in \cite{Jiang}.

\begin{lem}\Label{multifan}
Let $\tau$ be a cone in the multi-fan $\Delta_{\mathbf{\beta}}$.
If $\{\rho_{1},\cdots,\rho_{t}\}\subset link(\tau)$, and suppose
$\rho_{1},\cdots,\rho_{t}$ are contained in a cone
$\sigma\in\Delta_{\mathbf{\beta}}$. Then $\sigma\cup\tau$ is
contained in a cone of $\Delta_{\mathbf{\beta}}$.
\end{lem}

\begin{prop}\Label{moduleisom}
Let $\mathcal{M}(\mathcal{A})$ be the hypertoric DM stack associated to the stacky hyperplane arrangement $\mathcal{A}$, 
then we have an isomorphism of graded $A^{*}(\mathcal{M}(\mathcal{A}))$-modules:
$$\frac{\mathbb{Q}[\Delta_\mathbf{\beta}]}{Cir(\Delta_\mathbf{\beta})}\cong
\bigoplus_{(v,\sigma)\in Box(\Delta_{\mathbf{\beta}})}A^{*}(\mathcal{M}(\mathcal{A}(\sigma)))[deg(y^{(v,\sigma)})].$$
\end{prop}
\begin{pf}
We use arguments similar to those in Proposition 4.7 of
\cite{Jiang}. From Lemma \ref{lemma4.5} it is easy to see that 
$$\mathbb{Q}[\Delta_\mathbf{\beta}]\cong \bigoplus_{(v,\sigma)\in Box(\Delta_{\mathbf{\beta}})}y^{(v,\sigma)}\cdot \mathbb{Q}[M_{\beta}].$$
Consider $\bigoplus_{(v,\sigma)\in Box(\Delta_{\mathbf{\beta}})}y^{(v,\sigma)}\cdot
Cir(\Delta_\mathbf{\beta})$. It is an ideal of the ring $\bigoplus_{(v,\sigma)\in
Box(\Delta_{\mathbf{\beta}})}y^{(v,\sigma)}\cdot \mathbb{Q}[M_{\beta}]$,
so
$$\frac{\mathbb{Q}[\Delta_\mathbf{\beta}]}{Cir(\Delta_\mathbf{\beta})}\cong
\bigoplus_{(v,\sigma)\in Box(\Delta_{\mathbf{\beta}})}\frac{y^{(v,\sigma)}\cdot
\mathbb{Q}[M_{\beta}]}{y^{(v,\sigma)}\cdot Cir(\Delta_\mathbf{\beta})}.$$

For an element $(v,\sigma)\in Box(\Delta_{\mathbf{\beta}})$, let
$\rho_{1},\cdots,\rho_{l}\in link(\sigma)$. Then we
have an induced {\em matroid complex}
$M_{\beta(\sigma)}$, where $\beta(\sigma)$ is the map in the 
quotient stacky hyperplane arrangement $\mathcal{A}(\sigma)$ and the quotient extended stacky fan
$\mathbf{\Sigma/\sigma}$. Similarly from $\beta(\sigma)$, we have multi-fan 
$\Delta_{\beta(\sigma)}$  in $\overline{N(\sigma)}$. By Lemma
\ref{coarsecoh},
$A^{*}(\mathcal{M}(\mathcal{A}(\sigma)))\cong
\mathbb{Q}[M_{\beta(\sigma)}]/Cir(\Delta_{\beta(\sigma)})$.
For any element $(v,\sigma)\in
Box(\Delta_{\mathbf{\beta}})$, we construct an isomorphism
$$\Psi_{v}: \frac{\mathbb{Q}[M_{\beta(\sigma)}]}{Cir(\Delta_\mathbf{\beta(\sigma)})}[deg(y^{(v,\sigma)})]
\longrightarrow \frac{y^{(v,\sigma)}\cdot
\mathbb{Q}[M_{\beta}]}{y^{(v,\sigma)}\cdot Cir(\Delta_\mathbf{\beta})}.
$$
as follows. Consider the quotient  stacky hyperplane arrangement
$\mathcal{A}(\sigma)=(N(\sigma),\beta(\sigma),\theta(\sigma))$.
The hypertoric DM stack
$\mathcal{M}(\mathcal{A}(\sigma))$ is a closed
substack of $\mathcal{M}(\mathcal{A})$. Consider the morphism
$i:
\mathbb{Q}[y^{\widetilde{b}_{1}},\ldots,y^{\widetilde{b}_{l}}]\to
\mathbb{Q}[y^{b_{1}},\ldots,y^{b_{m}}]$ given by
$y^{\widetilde{b}_{i}}\mapsto y^{b_{i}}$. By Lemma
\ref{multifan}, it is easy to check that the ideal
$I_{M_{\beta(\sigma)}}$ is mapped to the ideal
$I_{M_{\beta}}$, so we have a morphism
$\mathbb{Q}[M_{\beta(\sigma)}]\to
\mathbb{Q}[M_{\beta}]$. Since $\mathbb{Q}[M_{\beta}]$ is a
subring of $\mathbb{Q}[\Delta_\mathbf{\beta}]$. Let
$\widetilde{\Psi}_{v}:
\mathbb{Q}[M_{\beta(\sigma)}][deg(y^{(v,\sigma)})]\to
y^{(v,\sigma)}\cdot \mathbb{Q}[M_{\beta}]$ be the morphism given by $y^{\widetilde{b}_{i}}\mapsto y^{(v,\sigma)}\cdot y^{b_{i}}$. Using
similar arguments as in Proposition 4.7 of \cite{Jiang}, we find that the ideal
$Cir(\Delta_\mathbf{\beta(\sigma)})$ goes to the ideal
$y^{(v,\sigma)}\cdot Cir(\Delta_\mathbf{\beta})$, so we have the morphism
$\Psi_{v}$ such that
$\Psi_{v}([y^{\widetilde{b}_{i}}])=[y^{(v,\sigma)}\cdot y^{b_{i}}]$.

Conversely, for $(v,\sigma)\in Box(\Delta_{\mathbf{\beta}})$, since
$\sigma$ is simplicial, for $\rho_{i}\subset \sigma$ we can choose $\theta_{i}\in
Hom_\mathbb{Z}(N,\mathbb{Q})$ such that $\theta_{i}(b_{i})=1$ and
$\theta_{i}(b_{i^{'}})=0$ for $b_{i^{'}}\neq b_{i}\in
\sigma$. We consider the following morphism $p:
\mathbb{Q}[y^{b_{1}},\ldots,y^{b_{m}}]\to
\mathbb{Q}[y^{\widetilde{b}_{1}},\ldots,y^{\widetilde{b}_{l}}]$ given by:
$$y^{b_{i}}\longmapsto\begin{cases}y^{\widetilde{b}_{i}}&\text{if $\rho_{i}\subseteq link(\sigma)$}\,,\\
-\sum_{j=1}^{l}\theta_{i}(b_{j})y^{\widetilde{b}_{j}}&\text{if
$\rho_{i}\subseteq \sigma$}\,,\\
0&\text{if $\rho_{i}\nsubseteq \sigma\cup
link(\sigma)$}\,.\end{cases}$$ Again by Lemma
\ref{multifan}, the ideal $I_{M_{\beta}}$ is mapped by $p$ to the ideal
$I_{M_{\beta(\sigma)}}$.  Then $p$ induces a
surjective map $\mathbb{Q}[M_{\beta}]\to
\mathbb{Q}[M_{\beta(\sigma)}]$ and a surjective
map $\widetilde{\Phi}_{v}: y^{(v,\sigma)}\cdot
\mathbb{Q}[M_{\beta}]\to
\mathbb{Q}[M_{\beta(\sigma)}][deg(y^{(v,\sigma)})]$. Using
the same computation as in Proposition 4.7 in \cite{Jiang}, the
relations $y^{(v,\sigma)}\cdot Cir(\Delta_\mathbf{\beta})$ is seen to go to
$Cir(\Delta_\mathbf{\beta(\sigma)})$. This yields another
morphism
$$\Phi_{v}: \frac{y^{(v,\sigma)}\cdot
\mathbb{Q}[M_{\beta}]}{y^{(v,\sigma)}\cdot
Cir(\Delta_\mathbf{\beta})}\longrightarrow
\frac{\mathbb{Q}[M_{\beta(\sigma)}]}{Cir(\Delta_\mathbf{\beta(\sigma)})}[deg(y^{(v,\sigma)})]
$$
so that $\Phi_{v}\Psi_{v}=1, \Psi_{v}\Phi_{v}=1$. So $\Psi_{v}$ is
an isomorphism.  We conclude by Lemma \ref{coarsecoh}.
\end{pf}

\subsection{The orbifold product}

In this section we compute the orbifold cup product. First for 
any $(v_1,\sigma_1),(v_2,\sigma_2)\in Box(\Delta_{\beta})$, we have the following lemma:

\begin{lem}\Label{triples}
If $\sigma_1\cup\sigma_2$ is a cone in the multi-fan $\Delta_{\beta}$, there 
exists a unique $(v_3,\sigma_3)\in Box(\Delta_{\beta})$ such that 
$\sigma_1\cup\sigma_2\cup\sigma_3$ is a cone in the multi-fan $\Delta_{\beta}$ and 
$v_1+v_2+v_3$ has no fractional part.
\end{lem}

\begin{pf}
Let $v_3=\lceil v_1+v_2\rceil_{\sigma_1\cup\sigma_2}- v_1-v_2$ and 
$\sigma_3$ the minimal cone in $\sigma_1\cup\sigma_2$ containing $v_3$. Then $(v_3,\sigma_3)$ satisfies the conditions of the Lemma. 
\end{pf}

The notation $(v_1,\sigma_1)+(v_2,\sigma_2)+(v_3,\sigma_3)\equiv 0$  
means that the triple $((v_1,\sigma_1),(v_2,\sigma_2),(v_3,\sigma_3))$ satisfies 
the conditions in Lemma \ref{triples}.

By \cite{CR2}, the 3-twisted sector
$\mathcal{M}(\mathcal{A})_{(g_1,g_2,g_3)}$ is the moduli space
of 3-pointed genus 0 degree 0 orbifold stable maps to
$\mathcal{M}(\mathcal{A})$. Let $\mathbb{P}^1(0, 1, \infty)$
be a genus 0 twisted curve with stacky structures possibly at $0,
1, \infty$. Consider a constant map 
$f: \mathbb{P}^1(0, 1, \infty) \to \mathcal{M}(\mathcal{A})$ with image $x\in \mathcal{M}(\mathcal{A})$.
This induces a morphism $$\rho: \pi_1^{orb}(\mathbb{P}^1(0, 1,
\infty))\to G_x, $$ where $G_x$ is the local group of the point
$x$. Let $\gamma_i$ be generators of $\pi_1^{orb}(\mathbb{P}^1(0,
1, \infty))$ and $g_i:=\rho(\gamma_i)$. The $g_i$ fixes the point
$x$. By Proposition \ref{fixbox}, $g_i$ corresponds to an element
$(v_i,\sigma_i)\in Box(\Delta_\mathbf{\beta})$.  
An argument similar to that in
Proposition 6.1 in \cite{BCS} shows that 3-twisted sectors of the
hypertoric DM stack $\mathcal{M}(\mathcal{A})$ are given by
\begin{equation}\Label{3-sector}
\coprod_{((v_{1},\sigma_{1}),(v_{2},\sigma_{2}),(v_{3},\sigma_{3}))\in Box(\Delta_{\mathbf{\beta}})^{3}, 
(v_{1},\sigma_1)+(v_{2},\sigma_2)+(v_{3},\sigma_3)\equiv 0}
\mathcal{M}(\mathcal{A}(\sigma_{123})),
\end{equation}
where $\sigma_{123}$ is the cone in $\Delta_{\mathbf{\beta}}$ satisfying 
$v_{1}+v_{2}+v_{3}=\sum_{\rho_{i}\subset
\sigma_{123}}a_{i}b_{i}$,
$a_{i}=1 , 2$.
Let $ev_{i}:
\mathcal{M}(\mathcal{A}(\sigma_{123}))\to
\mathcal{M}(\mathcal{A}(\sigma_{i}))$ be the
evaluation map. We have the obstruction bundle (see \cite{CR1}) $Ob_{(v_{1},v_{2},v_{3})}$ over the 3-twisted sector
$\mathcal{M}(\mathcal{A}(\sigma_{123}))$,
\begin{equation}\Label{obstruction}
Ob_{(v_{1},v_{2},v_{3})}=\left(e^{*}T\left(\mathcal{M}(\mathcal{A})\right)\otimes
H^{1}(C,\mathcal{O}_{C})\right)^{H}
\end{equation}
where $e:
\mathcal{M}(\mathcal{A}(\sigma_{123}))\to
\mathcal{M}(\mathcal{A})$ is the embedding,  $C\to
\mathbb{P}^{1}$ is the $H$-covering branching over three marked
points $\{0,1,\infty\}\subset \mathbb{P}^{1}$, and  $H$ is the
group generated by $v_{1},v_{2},v_{3}$.
Let $(v,\sigma)\in Box(\Delta_{\mathbf{\beta}})$, say $v\in N(\tau)$ for
some top dimensional cone $\tau$. Let $(\check{v},\sigma)\in
Box(\Delta_{\mathbf{\beta}})$ be the inverse of $v$ as an element
in the group $N(\tau)$. Equivalently, if
$\overline{v}=\sum_{\rho_{i}\subseteq \sigma}\alpha_{i}\overline{b}_{i}$
for $0<\alpha_{i}<1$, then $\check{\overline{v}}=\sum_{\rho_{i}\subseteq
\sigma}(1-\alpha_{i})\overline{b}_{i}$. 

\begin{lem}\label{equalto3}
Let $(v_{1},\sigma_{1}),(v_{2},\sigma_{2}),(v_{3},\sigma_{3})\in Box(\Delta_{\mathbf{\beta}})$ such that
$v_{1}+v_{2}+v_{3}\equiv 0$. Then if
$\overline{v}_{1}+\overline{v}_{2}+\overline{v}_{3}=\sum_{\rho_{i}\subseteq
\sigma_{123}}a_{i}\overline{b}_{i}$,
$\check{\overline{v}}_{1}+\check{\overline{v}}_{2}+\check{\overline{v}}_{3}=\sum_{\rho_{i}\subseteq
\sigma_{123}}c_{i}\overline{b}_{i}$,
where $a_{i},c_{i}=1~\text{or}~2$, then $a_{i}+c_{i}=2$ or $3$.
\end{lem}
\begin{pf}
Let  $\overline{v}_{i}=\sum_{\rho_{j}\subseteq \sigma_i} \alpha^{i}_j \overline{b}_j$, with $0<\alpha^{i}_j<1$ and
$i=1,2,3$. Then $\check{\overline{v}}_{i}=\sum_{\rho_{j}\subseteq \sigma_i} (1-\alpha^{i}_j) \overline{b}_j$. 
From the condition we have $\alpha^{1}_j+\alpha^{2}_j+\alpha^{3}_j=a_j=1~ \text{or} ~2$ and 
$(1-\alpha^{1}_j)+(1-\alpha^{2}_j)+(1-\alpha^{3}_j)=c_j=2~ \text{or} ~1$.
So if $\rho_j$ belongs to $\sigma_1, \sigma_2$ and $\sigma_3$, then 
$\alpha^{1}_j,\alpha^{2}_j,\alpha^{3}_j$ exist  and 
if $a_j=1~ \text{or} ~2$, then $c_j=2~ \text{or} ~1$. 
If  $\rho_j$ belongs to $\sigma_1, \sigma_2$, but not $\sigma_3$, then 
$\alpha^{3}_j$ doesn't exist and $\alpha^{1}_j+\alpha^{2}_j=a_j=1$, 
$(1-\alpha^{1}_j)+(1-\alpha^{2}_j)=c_j=1$. The cases that $\rho_j$ belongs to $\sigma_1,\sigma_3$
but  not $\sigma_2$, to $\sigma_2,\sigma_3$
but  not $\sigma_1$ are similar. We omit them.
\end{pf}

The stack $\mathcal{M}(\mathcal{A})$ is an abelian  DM stack,
i.e. the local groups are all abelian groups. For any
3-twisted sector
$\mathcal{M}(\mathcal{A}(\sigma_{123}))$,
the normal bundle
$N(\mathcal{M}(\mathcal{A}(\sigma_{123}))/\mathcal{M}(\mathcal{A}))$
can be split into the direct sum of some line bundles under the
group action. It follows from the definition that if
$\overline{v}=\sum_{\rho_{i}\subseteq\sigma_{123}}\alpha_{i}\overline{b}_{i}$,
then the action of $v$ on the normal bundle
$N(\mathcal{M}(\mathcal{A}(\sigma_{123}))/\mathcal{M}(\mathcal{A}))$
is given by the diagonal matrix $diag(\alpha_{i},1-\alpha_{i})$. A
general result in \cite{CH} and \cite{JKK} about the obstruction
bundle and Lemma \ref{equalto3} imply the following Proposition.

\begin{prop}\label{obstructionbdle}
Let $\mathcal{M}(\mathcal{A})_{(v_{1},v_{2},v_{3})}=
\mathcal{M}(\mathcal{A}(\sigma_{123}))$
be a 3-twisted sector of the stack $\mathcal{M}(\mathcal{A})$ such that 
$v_{1},v_{2},v_{3}\neq 0$.
Then the Euler class of the obstruction bundle
$Ob_{(v_{1},v_{2},v_{3})}$ on
$\mathcal{M}(\mathcal{A}(\sigma_{123}))$
is
$$\prod_{a_{i}=2}c_{1}(L_{i})|_{\mathcal{M}(\mathcal{A}(\sigma_{123}))}
\cdot\prod_{a_{i}=1,\,\, \alpha^{1}_j,\alpha^{2}_j,\alpha^{3}_j\text{ exist}}c_{1}(L^{-1}_{i})|_{\mathcal{M}(\mathcal{A}(\sigma_{123}))},$$
where $L_{i}$ is the line bundle over
$\mathcal{M}(\mathcal{A})$ defined in Definition \ref{linebdle}. $\square$
\end{prop}

To prove the main theorem, we introduce two Lemmas.  
For any two pairs $(c_{1},\sigma_{1}), (c_{2},\sigma_{2})\in N^{\Delta_\mathbf{\beta}}$,
there exist two unique elements $(v_{1},\tau_{1}), (v_{2},\tau_{2})\in Box(\Delta_\mathbf{\beta})$ such that 
$\tau_1\subseteq\sigma_1, \tau_2\subseteq\sigma_2$ and $c_{1}=v_{1}+\sum_{\rho_{i}\subseteq \sigma_{1}}m_{i}b_{i}$, 
$c_{2}=v_{2}+\sum_{\rho_{i}\subseteq \sigma_{2}}n_{i}b_{i}$, where 
$m_{i},n_{i}$ are nonnegative integers. Let $(v_{3},\sigma_{3})$ be the unique element in $Box(\Delta_{\mathbf{\beta}})$  such that 
$v_{1}+v_{2}+v_{3}\equiv 0$ in the local group given by $\sigma_{1}\cup\sigma_{2}$. 
Let $\overline{v}_{i}=\sum_{\rho_{j}\subseteq \sigma_i} \alpha^{i}_j \overline{b}_j$, with $0<\alpha^{i}_j<1$ and
$i=1,2,3$. The existence of  $\alpha^{1}_j,\alpha^{2}_j,\alpha^{3}_j$  means that $\rho_j$ belongs to $\sigma_1, \sigma_2, \sigma_3$.
Let $\sigma_{123}$ be the cone in $\Delta_{\mathbf{\beta}}$ such that 
$\overline{v}_{1}+\overline{v}_{2}+\overline{v}_{3}=\sum_{\rho_i\subseteq\sigma_{123}}a_{i}\overline{b}_i$, with $a_i=1$ or $2$. Let $I$ be the set of 
$i$ such that $a_i=1$ and  $\alpha^{1}_j, \alpha^{2}_j, \alpha^{3}_j$ exist, $J$ the set of $j$ such that $\rho_j$ belongs to $\sigma_{123}$ but not to $\sigma_3$.
We have the following lemma for the ceiling functions:

\begin{lem}\Label{Gauss}
$\lceil c_1\rceil_{\sigma_1}+\lceil c_2\rceil_{\sigma_2}-\lceil c_1+c_2 \rceil_{\sigma_1\cup\sigma_2}=
\lceil v_1\rceil_{\tau_1}+\lceil v_2\rceil_{\tau_2}-\lceil v_1+v_2 \rceil_{\tau_1\cup\tau_2}$.
\end{lem}

\begin{pf}
By the definition of ceiling functions, we have $\lceil c_{1}\rceil_{\sigma_1}=\lceil v_{1}\rceil_{\tau_1}+\sum_{\rho_{i}\subseteq \sigma_{1}}m_{i}b_{i}$ and $\lceil c_{2}\rceil_{\sigma_2}=\lceil v_{2}\rceil_{\tau_2}+\sum_{\rho_{i}\subseteq \sigma_{2}}n_{i}b_{i}$. The Lemma follows.
\end{pf}

\begin{lem}\Label{keylem}
If $\sigma_{1}\cup\sigma_{2}$ is a cone in $\Delta_\mathbf{\beta}$ for the two pairs
$(c_{1},\sigma_{1}), (c_{2},\sigma_{2})$, then the product $y^{(c_1,\sigma_1)}\cdot y^{(c_2,\sigma_2)}$ in 
(\ref{product}) can be given by
\begin{equation}\Label{productmain}
\begin{cases}
(-1)^{|I|+|J|}y^{(c_{1}+c_{2}+\sum_{i\in I}b_{i}+\sum_{i\in J}b_{i},\sigma_{1}\cup\sigma_{2})}&\text{if
$\overline{v}_{1},\overline{v}_{2}\neq 0$ and $\overline{v}_{1}\neq\check{\overline{v}}_{2}$}\,,\\
(-1)^{|J|}y^{(c_{1}+c_{2}+\sum_{i\in J}b_{i},\sigma_{1}\cup\sigma_{2})}&\text{if
$\overline{v}_{1},\overline{v}_{2}\neq 0$ and $\overline{v}_{1}=\check{\overline{v}}_{2}$}\,,\\
y^{(c_{1}+c_{2},\sigma_{1}\cup\sigma_{2})} &\text{if
$\overline{v}_{1}$ or $\overline{v}_{2}= 0$}\,.
\end{cases} 
\end{equation}
\end{lem}

\begin{pf}
First for a fixed ray $\rho_i$ and $0<\alpha_{1},\alpha_{2}<1$, by the definition of ceiling functions, we find that 
\begin{equation}\Label{check}
\lceil \alpha_{1}b_{i}\rceil_{\rho_i}+\lceil \alpha_{2}b_{i}\rceil_{\rho_i}-\lceil \alpha_{1}b_{i}+\alpha_{2}b_{i}\rceil_{\rho_i}=
\begin{cases}
0&\text{if
$\alpha_1+\alpha_2>1$}\,,\\
b_{i}
&\text{if
$\alpha_1+\alpha_2\leq 1$}\,.
\end{cases} 
\end{equation}
Since $\epsilon(c_1,c_2)=\lceil c_1\rceil_{\sigma_1}+\lceil c_2\rceil_{\sigma_2}-\lceil c_1+c_2\rceil_{\sigma_1\cup\sigma_2}$, 
by Lemma \ref{Gauss}, we need to check that 
$$
\lceil v_1\rceil_{\tau_1}+\lceil v_2\rceil_{\tau_2}-\lceil v_1+v_2\rceil_{\tau_1\cup\tau_2}=
\begin{cases}
\sum_{i\in I}b_{i}+\sum_{i\in J}b_{i}&\text{if
$\overline{v}_{1},\overline{v}_{2}\neq 0$ and $\overline{v}_{1}\neq\check{\overline{v}}_{2}$}\,,\\
\sum_{i\in J}b_{i}&\text{if
$\overline{v}_{1},\overline{v}_{2}\neq 0$ and $\overline{v}_{1}=\check{\overline{v}}_{2}$}\,,\\
0 &\text{if
$\overline{v}_{1}$ or $\overline{v}_{2}= 0$}\,.
\end{cases}
$$ 
This can be easily proven using (\ref{check})  and Lemma \ref{equalto3}.
\end{pf}

\subsection{Proof of Theorem 1.1}
By Proposition \ref{moduleisom}, it remains to prove that the
orbifold cup product is the same as the product in the  ring
$\mathbb{Q}[\Delta_\mathbf{\beta}]$. By Lemma \ref{keylem}, we need to prove that the
orbifold cup product is the same as the product in (\ref{productmain}).  It suffices to consider the
canonical generators $y^{b_{i}}$, $y^{(v,\sigma)}$ for $(v,\sigma)\in
Box(\Delta_{\mathbf{\beta}})$.

Consider $y^{(v,\sigma)}\cup_{orb}y^{b_{i}}$ with $(v,\sigma)\in
Box(\Delta_{\mathbf{\beta}})$. The element $(v,\sigma)$ determines a
twisted sector
$\mathcal{M}(\mathcal{A}(\sigma))$. The
corresponding twisted sector to $b_{i}$ is the whole hypertoric
stack  $\mathcal{M}(\mathcal{A})$. It is easy to see that the
3-twisted sector relevant to this product is
$\mathcal{M}(\mathcal{A})_{(v,1,v^{-1})}\cong
\mathcal{M}(\mathcal{A}(\sigma))$, where $v^{-1}$
denotes the inverse of $v$ in the local group. It follows from the
dimension formula in \cite{CR1} that the obstruction bundle over
$\mathcal{M}(\mathcal{A})_{(v,1,v^{-1})}$ has rank zero. It is
immediate from definition that
$y^{(v,\sigma)}\cup_{orb}y^{b_{i}}=y^{(v+b_{i},\sigma\cup \rho_{i})}$ if 
there is a cone in $\Delta_{\beta}$ containing $\overline{v},\overline{b}_{i}$.
This is the third case  in (\ref{productmain}).

Now consider $y^{(v_{1},\sigma_{1})}\cup_{orb}y^{(v_{2},\sigma_{2})}$, where $(v_{1},\sigma_{1}),(v_{2},\sigma_{2})\in
Box(\Delta_{\mathbf{\beta}})$. By (\ref{3-sector}), we see that
if $\sigma_1\cup\sigma_2$ is not a cone in $\Delta_{\mathbf{\beta}}$,  then there is no 3-twisted
sector corresponding to the elements $v_{1},v_{2}$. Thus the
product is zero by definition. On the other hand, by definition of the ring $\mathbb{Q}[\Delta_\mathbf{\beta}]$, we have $y^{(v_{1},\sigma_{1})}\cdot y^{(v_{2},\sigma_{2})}=0$. So
$y^{(v_{1},\sigma_{1})}\cup_{orb}y^{(v_{2},\sigma_{2})}=y^{(v_{1},\sigma_{1})}\cdot y^{(v_{2},\sigma_{2})}$. If 
$\sigma_1\cup\sigma_2$ is a cone in $\Delta_{\mathbf{\beta}}$, let $(v_{3},\sigma_{3})\in
Box(\Delta_{\mathbf{\beta}})$ such that $\overline{v}_{3}\in
\sigma_{123}$ and $v_{1}v_{2}v_{3}=1$
in the local group. Then we have the 3-twisted sector
$\mathcal{M}(\mathcal{A}(\sigma_{123}))$.
Let $ev_{i}:
\mathcal{M}(\mathcal{A}(\sigma_{123}))\to
\mathcal{M}(\mathcal{A}(\sigma_{i}))$ be the
evaluation maps. The element $y^{(v,\sigma)}$ is the class $1$ in the
cohomology of the twisted sector
$\mathcal{M}(\mathcal{A}(\sigma))$. From the
definition of orbifold cup product \cite{CR1}, \cite{AGV}, we have:
$$y^{(v_{1},\sigma_{1})}\cup_{orb}y^{(v_{2},\sigma_{2})}=(\breve{ev}_{3})_{*}(ev_{1}^{*}y^{(v_{1},\sigma_{1})}
\cdot ev_{2}^{*}y^{(v_{2},\sigma_{2})}\cdot e(Ob_{(v_{1},v_{2},v_{3})})),$$
where $\breve{ev}_{3}=\mathcal{I}\circ ev_{3}:
\mathcal{M}(\mathcal{A}(\sigma_{123}))\to
\mathcal{M}(\mathcal{A})_{(\check{v}_{3})}$ is the composite of
$ev_{3}$ and the natural involution $\mathcal{I}:
\mathcal{M}(\mathcal{A})_{(v_{3})}\to
\mathcal{M}(\mathcal{A})_{(\check{v}_{3})}$. Let 
$\overline{v}_{i}=\sum_{\rho_{j}\subseteq \sigma_i} \alpha^{i}_j \overline{b}_j$, with $0<\alpha^{i}_j<1$ and
$i=1,2,3$. 
Let $I$ denote the set of $i$
such that $a_{i}=1$ and $\alpha^{1}_j,\alpha^{2}_j,\alpha^{3}_j$ exist,   $J$  the set of $j$
such that $\rho_{i}$ belongs to
$\sigma_{123}$, but not belong to
$\sigma_{3}$.  

If some $\overline{v}_{i}=0$, for example, $\overline{v}_{1}=0$, then $v_{1}$ is a torsion 
element in $N$ which means that the action of $v_{1}$ is trivial on the hypertoric DM stack.
Then the 3-twisted sector corresponding to $v_{1},v_{2}$
is isomorphic to the twisted sector $\mathcal{M}(\mathcal{A}(\sigma_{2}))$ and the obstruction 
bundle over $\mathcal{M}(\mathcal{A}(\sigma_{2}))$ is zero by \cite{CR1}. In this case 
the set $I$ and $J$ are all empty. So $y^{(v,\sigma_{1})}\cup_{orb}y^{(v,\sigma_{2})}=y^{(v_{1}+v_{2},\sigma_{1}\cup \sigma_{2})}$.
This is again the third case  in (\ref{productmain}).

Now we assume that $\overline{v}_{1},\overline{v}_{2}\neq 0$.
If $\overline{v}_{1}=\check{\overline{v}}_{2}$, then $\overline{v}_{3}=0$, $\sigma_{123}=\sigma_1$ and $v_1+v_2=\sum_{\rho_j\subseteq \sigma_1}b_{j}$. 
So the 3-twisted sector corresponding to $v_{1},v_{2}$
is isomorphic to the twisted sector $\mathcal{M}(\mathcal{A}(\sigma_{1}))$ and the obstruction 
bundle over $\mathcal{M}(\mathcal{A}(\sigma_{1}))$ is zero by \cite{CR1}. The set $J$ is the set 
$j$ such that $\rho_j\subseteq \sigma_1$. So we have
\begin{eqnarray}
\nonumber y^{(v_{1},\sigma_{1})}\cup_{orb}y^{(v_{2},\sigma_{2})}&=&
y^{0}\cdot \prod_{i\in
J}y^{b_{i}}\cdot \prod_{i\in
J}(-y^{b_{i}})  \\ \nonumber
&=&
(-1)^{|J|}\cdot y^{(v_{1}+v_2+\sum_{i\in J}b_{i},
\sigma_{1}\cup\sigma_{2})}, \nonumber
\end{eqnarray}
which is the second case in (\ref{productmain}). 

If $\overline{v}_{1}\neq \check{\overline{v}}_{2}$, then
$\overline{v}_{3}\neq 0$ and the obstruction bundle over the 3-twisted sector $\mathcal{M}(\mathcal{A}(\sigma_{123}))$
is given by Proposition 
\ref{obstructionbdle}. So we have:
$$y^{(v_{1},\sigma_{1})}\cup_{orb}y^{(v_{2},\sigma_{2})}=y^{(\check{v}_{3},\sigma_{3})}
\cdot\prod_{a_{i}=2}y^{b_{i}}\cdot\prod_{i\in
J}y^{b_{i}}\cdot \prod_{i\in I}(-y^{b_{i}})\cdot \prod_{i\in
J}(-y^{b_{i}}).$$ 
Since
$\check{v}_{3}+\sum_{a_{i}=2}b_{i}+\sum_{i\in J}b_{i}=v_{1}+v_{2}$,  
we have
\begin{eqnarray}
\nonumber y^{(v_{1},\sigma_{1})}\cup_{orb}y^{(v_{2},\sigma_{2})}&=&(-1)^{|I|+|J|}\cdot y^{(v_{1}+v_2,\sigma_{1}\cup\sigma_{2})}
\cdot \prod_{i\in I}y^{b_{i}}\cdot \prod_{i\in
J}y^{b_{i}} \\ \nonumber
&=&
(-1)^{|I|+|J|}\cdot y^{(v_{1}+v_2+\sum_{i\in I}b_{i}+\sum_{i\in
J}b_{i},\sigma_{1}\cup\sigma_{2})}, \nonumber
\end{eqnarray}
which is the first case in (\ref{productmain}).
~$\square$

\section{Applications}\label{app}
In this section we compute some examples of the orbifold Chow
rings of  hypertoric DM stacks. In particular, we relate the
hypertoric stack to crepant resolutions.

Let $N=\mathbb{Z}$ and $\Sigma$ the fan of projective line
$\mathbb{P}^{1}$ generated by $\{(1),(-1)\}$. Let $\beta:
\mathbb{Z}^{n}\to N$ be the map given by
$b_{1}=(1),b_{2}=(-1)$ and $b_{i}=(1)$ for $i\geq 2$. Consider the
following exact sequences $$0\longrightarrow
\mathbb{Z}^{n-1}\longrightarrow
\mathbb{Z}^{n}\stackrel{\beta}{\longrightarrow} N\longrightarrow
0\longrightarrow 0,$$
$$0\longrightarrow \mathbb{Z}\longrightarrow
\mathbb{Z}^{n}\stackrel{\beta^{\vee}}{\longrightarrow}
\mathbb{Z}^{n-1}\longrightarrow 0\longrightarrow 0,$$ where the
Gale dual $\beta^{\vee}$ is given by the column vectors of the  matrix $$A=\left[
\begin{array}{cccccc}
  1&1&0&0&\cdots&0 \\
  1&0&-1&0&\cdots&0 \\
  1&0&0&-1&\cdots&0 \\
  \vdots&\vdots&\vdots&\vdots&\ddots&\vdots \\
  1&0&0&0&\cdots&-1\\
\end{array}
\right].$$ Note that $A$ is unimodular in the sense of 
\cite{HS}. Taking $Hom_\mathbb{Z}(-,\mathbb{C}^{\times})$ yields
$$1\longrightarrow
(\mathbb{C}^{\times})^{n-1}\stackrel{\alpha}{\longrightarrow}
(\mathbb{C}^{\times})^{n}\longrightarrow
\mathbb{C}^{\times}\longrightarrow 1.$$ So
$G=(\mathbb{C}^{\times})^{n-1}$.  Choose 
$\theta=(1,1,\cdots,1)$  in
$\mathbb{Z}^{n-1}$, then it is a generic element. The extended stacky fan 
$\mathbf{\Sigma}=(N,\Sigma,\beta)$ is induced from the stacky hyperplane 
arrangement $\mathcal{A}=(N,\beta,\theta)$, where $\mathcal{H}$ is the 
hyperplane arrangement whose normal fan is $\Sigma$. It
is easy to see that the  toric DM stack is the projective
line $\mathbb{P}^{1}$. The hypertoric DM stack is the crepant
resolution of the Gorenstein orbifold
$[\mathbb{C}^{2}/\mathbb{Z}_{n}]$. To see this, from the construction of
hypertoric DM stack, we have:
\begin{equation}\label{resolution}
1\longrightarrow
(\mathbb{C}^{\times})^{n-1}\stackrel{\alpha^{L}}{\longrightarrow}
(\mathbb{C}^{\times})^{2n}\longrightarrow
(\mathbb{C}^{\times})^{n+1}\longrightarrow 1,\end{equation} where
$\alpha^{L}$ is given by the matrix $[\beta^{\vee},-\beta^{\vee}]$. Let
$\mathbb{C}[z_1,...,z_n,w_1,...,w_n]$ be the coordinate ring of
$\mathbb{C}^{2n}$. So the ideal $I_{\beta^{\vee}}$ in (\ref{ideal1}) is
generated by the following equations:
$$\begin{cases}z_{1}w_{1}+z_{2}w_{2}=0\,,\\
z_{1}w_{1}-z_{3}w_{3}=0\,,\\
\cdots \cdots\cdots\cdots\cdots\\
z_{1}w_{1}-z_{n}w_{n}=0\,.\end{cases}$$ 
Hence $Y$ is the subvariety of
$\mathbb{C}^{2n}-V(\mathcal{I}_{\theta})$ determined by the above
ideal. The action of $G$ on $Y$ is through the map $\alpha^{L}$ in
(\ref{resolution}). The hypertoric DM stack associated to
$\mathcal{A}$ is  $\mathcal{M}(\mathcal{A})=[Y/G]$. From
Proposition \ref{orientation}, the hypertoric DM stack is
independent to the coorientations of the hyperplanes. This means
that we can give the stacky hyperplane arrangement  $\mathcal{A}$ as
follows. Let $b_{i}=1$ for $1\leq i\leq n$.
Then the Gale dual map $\beta^{\vee}:
\mathbb{Z}^{n}\to \mathbb{Z}^{n-1}$ is given by the
matrix
$$A=\left[
\begin{array}{cccccc}
  1&-1&0&0&\cdots&0 \\
  0&1&-1&0&\cdots&0 \\
  0&0&1&-1&\cdots&0 \\
  \vdots&\vdots&\ddots&\ddots&\ddots&\vdots \\
  0&0&0&\cdots&1&-1\\
\end{array}
\right],$$ which is exactly the matrix in Lemma 10.2 in \cite{HS},
from which it follows that the coarse moduli space $Y(\beta^{\vee},\theta)$ of
$\mathcal{M}(\mathcal{A})=[Y/G]$ is the crepant resolution of
the Gorenstein orbifold $[\mathbb{C}^{2}/\mathbb{Z}_{n}]$. The
core of the hypertoric DM stack $\mathcal{M}(\mathcal{A})$ is
a chain of $n-1$ copies of $\mathbb{P}^{1}$ with normal crossing
divisors corresponding to the multi-fan
$\Delta_{\mathbf{\beta}}$.

\begin{rmk}
This is an example of \cite{Kro}, in which it is shown that the
minimal resolution of a surface singularity of ADE type can be
constructed as a hyperk\"ahler quotient.
\end{rmk}

The $\mathbb{Z}_n$-action defining the Gorenstein orbifold $[\mathbb{C}^{2}/\mathbb{Z}_{n}]$ is given by $\lambda(x,y)=(\lambda x,\lambda^{-1}y)$ for $\lambda\in \mathbb{Z}_n$. There are $n-1$
twisted sectors each of which is isomorphic to
$\mathcal{B}\mathbb{Z}_{n}$ with age $1$. There are only
dimension zero cohomology on the untwisted sector and twisted
sectors. So we prove the following Proposition:

\begin{prop}\label{orbifold}
The orbifold Chow ring
$A^{*}_{orb}([\mathbb{C}^{2}/\mathbb{Z}_{n}])$ of
$[\mathbb{C}^{2}/\mathbb{Z}_{n}]$ is isomorphic to the ring
$$\frac{\mathbb{C}[x_{1},\cdots,x_{n-1}]}{\{x_{i}x_{j}: 1\leq i,j\leq
n-1\}}.$$  
\end{prop}
Since the crepant resolution is a manifold, the orbifold Chow
ring is the ordinary Chow ring. By Theorem 1.1, we have
\begin{prop}\label{crepant}
The  Chow ring of $\mathcal{M}(\mathcal{A})$ is isomorphic to
the ring
$$\frac{\mathbb{C}[y_{1},\cdots,y_{n-1}]}{\{y_{i}y_{j}: 1\leq i,j\leq n-1\}},$$
which is isomorphic to the orbifold cohomology ring of the
Gorenstein orbifold $[\mathbb{C}^{2}/\mathbb{Z}_{n}]$.
\end{prop}

\begin{pf}
By Theorem 1.1, the Chow ring of $\mathcal{M}(\mathcal{A})$ is
isomorphic to the ring:
$$\frac{\mathbb{C}[y_{1},\cdots,y_{n}]}{\{y_{1}-y_{n}+y_{3}+\cdots+y_{n-1},y_{i}y_{j}: 1\leq i,j\leq n-1\}}$$ which we can
easily check that  this ring is isomorphic to the orbifold
cohomology ring of $[\mathbb{C}^{2}/\mathbb{Z}_{n}]$ in
Proposition \ref{orbifold}.
\end{pf}

Y. Ruan \cite{R} conjectured that, among other things, the orbifold cohomology ring of
a hyperk\"ahler orbifold is isomorphic to the ordinary cohomology
ring of a hyperk\"ahler resolution (which is crepant).   For the
orbifold $[\mathbb{C}^{2}/\mathbb{Z}_{n}]$, the crepant resolution
$Y(\beta^{\vee},\theta)$ is smooth, we have that
$\mathcal{M}(\mathcal{A})\cong Y(\beta^{\vee},\theta)$. From Proposition
\ref{crepant}, the conjecture is true.  

A conjecture equating Gromov-Witten theories of an orbifold and its crepant resolutions, as proposed in \cite{BG}, is recently proven in genus 0 for $[\mathbb{C}^2/\mathbb{Z}_3]$, see \cite{BGP}. The comparison of two Gromov-Witten theories requires certain change of variables. For $[\mathbb{C}^2/\mathbb{Z}_3]$ case, see \cite{BGP}. For $[\mathbb{C}^2/\mathbb{Z}_4]$ case the following change of variables is found in \cite{BJ}:
$$\begin{cases}y_{1}=\frac{1}{4}(\sqrt{2}x_{1}+2ix_{2}-\sqrt{2}x_{3})\,,\\
y_{2}=\frac{1}{4}(\sqrt{2}ix_{1}-2ix_{2}+\sqrt{2}ix_{3})\,,\\
y_{3}=\frac{1}{4}(-\sqrt{2}x_{1}+2ix_{2}+\sqrt{2}x_{3})\,.\end{cases}$$
Under this change of variables, the genus zero Gromov-Witten
potential of the crepant resolution is seen to coincide with  the
genus zero orbifold Gromov-Witten potential of the orbifold
$[\mathbb{C}^{2}/\mathbb{Z}_{4}]$, see \cite{BJ}.

For a toric orbifold, it is known that adding rays in the simplicial fan can give a crepant resolution. In the end of the paper we compute an example and explain that adding rays in the stacky hyperplane arrangement doesn't give a smooth hypertoric variety
which means that it is not easy in general to give a crepant resolution in hyperk\"ahler geometry. 

\begin{example}
Let $\mathbf{\Sigma}=(N,\Sigma,\beta)$ be an extended stacky fan,
where $N=\mathbb{Z}^{2}$, the simplicial
fan $\Sigma$ is the fan of weighted projective plane
$\mathbb{P}(1,1,2)$, and $\beta: \mathbb{Z}^{3}\to N$
is given by the vectors
$\{b_{1}=(1,0),b_{2}=(0,1),b_{3}=(-1,-2),\}$, where
$b_{1},b_{2},b_{3}$ are the generators of the rays in $\Sigma$.  
The  generic element $\theta=(1)\in DG(\beta)\cong
\mathbb{Z}$ determines the fan $\Sigma$.  The stacky hyperplane arrangement 
$\mathcal{A}=(N,\beta,\theta)$ induces $\mathbf{\Sigma}$. The hypertoric DM stack
is $\mathcal{M}(\mathcal{A})=T^*(\mathbb{P}(1,1,2))$. 
From Theorem 1.1,
$$A_{orb}^{*}(\mathcal{M}(\mathcal{A}))\cong 
\frac{\mathbb{Q}[x_{1},x_{2},x_{3},x_{4}]}{(x_{1}-x_{3},x_{2}-2x_{3},x_{4}^{2}, x_{1}x_{2}x_{3},x_{4}x_{2},x_{4}x_{1}x_{3})}
\cong \frac{\mathbb{Q}[x_{3},x_{4}]}{(x_{4}^{2},x_{3}^{3},x_{3}x_{4})}.$$

Let $b_{4}=(0,-1)$ and consider the new map $\beta': \mathbb{Z}^{4}\to N$
which is given by the vectors
$\{b_{1},b_{2},b_{3},b_{4}\}$. Choose generic element $\theta'=(1,1)\in \mathbb{Z}^{2}=DG(\beta')$ and we get a new 
stacky hyperplane arrangement $\mathcal{A}'=(N,\beta',\theta')$ which induces the extended stacky fan  $\mathbf{\Sigma'}=(N,\Sigma,\beta')$.
The hypertoric DM stack 
$\mathcal{M}(\mathcal{A}')$ is the stack corresponding to  $\mathcal{A}'$.
From the definition of Box, $(\frac{1}{2}b_{1}+\frac{1}{2}b_3, \rho_1+\rho_3)$ is again a box element which 
determines a twisted sector. We compute that $A^{*}_{orb}(\mathcal{M}(\mathcal{A}'))$ is isomorphic to 
$$ 
\frac{\mathbb{Q}[x_{1},x_{2},x_{3},x_{4},v]}{(x_{1}-x_{3},x_{2}-2x_{3}-x_{4},x_{2}x_{4},x_{1}x_{2}x_{3},x_{1}x_{3}x_{4},v^{2},vx_2,vx_4)}
\cong \frac{\mathbb{Q}[x_{3},x_{4},v]}{(x_{3}x_{4}+x_{4}^{2},x_{3}^{3},x_{3}^{2}x_{4},v^{2},vx_3,vx_4)}.$$
We check that  $A_{orb}^{*}(\mathcal{M}(\mathcal{A}))$ is not isomorphic to 
the ring $A^{*}_{orb}(\mathcal{M}(\mathcal{A}'))$. 

We give two comments here. First, if the crepant resolution conjecture 
is true, then  $\mathcal{M}(\mathcal{A}')$ is not a hyperk\"ahler resolution.
On the other hand,
the map $\beta$ is given by the matrix $B=\left[
\begin{array}{ccc}
  1&0&-1 \\
  0&1&-2 
\end{array}
\right]$  and the map $\beta^{'}$ given by $B^{'}=\left[
\begin{array}{cccc}
  1&0&-1&0 \\
  0&1&-2&-1 
\end{array}
\right]$ It is easy to see that $B^{'}$ is not unimodular which means that $\mathcal{M}(\mathcal{A}')$ is 
not smooth, see \cite{HS}. So adding columns in $B$ can't make it unimodular which means that adding rays in the stacky hyperplane arrangement can not give a 
hyperk\"ahler resolution. 
Second, since $\mathcal{M}(\mathcal{A}')$ is still a hypertoric orbifold, we wish that $\mathcal{M}(\mathcal{A}')$
is a partial resolution and there is a morphism $\mathcal{M}(\mathcal{A}^{'})\longrightarrow \mathcal{M}(\mathcal{A})$.
But this is impossible because the two rings have different dimensions and it would violate the McKay correspondence
statement if there exists a morphism $\mathcal{M}(\mathcal{A}^{'})\longrightarrow \mathcal{M}(\mathcal{A})$, see 
\cite{Ya1},\cite{Ya2}.
\end{example}

$\mathbf{Question:}$ Is there a combinatorial description of a hyperk\"ahler resolution of hypertoric orbifolds?


\end{document}